\documentclass[10pt]{article}

\usepackage{amsmath}
\usepackage{amssymb}
\usepackage{amsthm}

\usepackage{etex}
\usepackage{bm}
\usepackage{natbib}
\usepackage[plain,noend]{algorithm2e}
\usepackage{subcaption}
\usepackage{placeins}
\usepackage{verbatim}
\usepackage{amsbsy}
\usepackage{booktabs}
\usepackage{array}
\usepackage{multirow}
\usepackage{floatrow}
\usepackage{tabularx}
\usepackage[pdftex]{graphicx}
\usepackage{tabularx,ragged2e,caption}
\usepackage{indentfirst}
\usepackage{float}
\usepackage{enumerate}
\usepackage{pgf}
\usepackage{tikz}
\usepackage[utf8]{inputenc}

\usetikzlibrary{arrows,automata}
\usetikzlibrary{positioning,backgrounds}
\usepackage{marginnote}

\newtheorem{thm}{Theorem}
\newtheorem{defi}{Definition}

\newcommand{\altmean}{\eta}
\newcommand{\altcov}{R}
\newcommand{\ept}{\xi}

\newcommand{\para}{\theta}
\newcommand{\psup}{\Theta}
\newcommand{\psupalt}{{\Theta_1}}
\newcommand{\psupnull}{{\Theta_0}}
\newcommand{\psupi}{{\Theta_i}}
\newcommand{\mlealt}{^{\text{\tiny $H_1$}}_{\text{\tiny MLE}}}
\newcommand{\mlenull}{^{\text{\tiny $H_0$}}_{\text{\tiny MLE}}}
\newcommand{\mlei}{^{\text{\tiny $H_i$}}_{\text{\tiny MLE}}}
\newcommand{\nuismlealt}{\nuis_{1,\text{\tiny MLE}}}
\newcommand{\nuismlenull}{\nuis_{0,\text{\tiny MLE}}}
\newcommand{\nuismlei}{\nuis_{i,\text{\tiny MLE}}}
\newcommand{\info}{\mathcal{I}}

\newcommand{\dinfo}{\mathcal{I}^T_\mathcal{V}}

\newcommand{\fracinfo}{\mathcal{FI}_\mathcal{V}^T}

\newcommand{\decinfo}{\mathcal{I}_\mathcal{V}^T}

\newcommand{\bdecinfo}{\mathcal{I}^{T}_\mathcal{V}}

\newcommand{\bdinfofull}{\mathcal{I}^{T}_\mathcal{V}}

\newcommand{\bdinfo}{\mathcal{I}^{T}_\mathcal{V}}
\newcommand{\uncp}{U}

\newcommand{\cn}{C_\declr} 
\newcommand{\declr}{\mathcal{V}} 

\newcommand{\decz}{\mathcal{V}(1)}   
\newcommand{\dec}{W}

\newcommand{\prior}{\pi}
\newcommand{\bfac}{\text{\normalfont BF}}
\newcommand{\zlr}{z}

\newcommand{\post}{p}
\newcommand{\ranx}{X}
\newcommand{\datx}{x}
\newcommand{\xsup}{\mathcal{X}}

\newcommand{\tk}{TK-optimality}
\newcommand{\tkc}{TK-optimality criterion}
\newcommand{\tkcs}{TK-optimality criteria}
\newcommand{\tkw}{TK-optimal}
\newcommand{\popt}{P-optimality}
\newcommand{\poptc}{P-optimality criterion}

\newcommand{\poptw}{P-optimal}
\newcommand{\dopt}{D-optimality}

\newcommand{\doptw}{D-optimal}

\newcommand{\crit}{\phi}

\newcommand{\critpb}{\phi_{\postbayes}}
\newcommand{\critpl}{\phi_{\postlrt}}
\newcommand{\crittb}{\phi_{\tkbayes}}
\newcommand{\crittl}{\phi_{\tklrt}}

\newcommand{\covoff}{\sigma_{\text{is}}}

\newcommand{\lr}{\text{\normalfont LR}}

\newcommand{\ri}{\mathcal{R}I}

\newcommand{\ob}{_{\text{\normalfont obs}}}
\newcommand{\co}{_{\text{\normalfont com}}}
\newcommand{\mis}{_{\text{\normalfont mis}}}
\newcommand{\var}{\text{Var}}

\newcommand{\Epara}{E_{\para}}

\newcommand{\Ex}{E_X}

\newcommand{\Exone}{E_{X_1}}
\newcommand{\Extwo}{E_{X_2}}

\newcommand{\Exmis}{E_{\ranx\mis}}

\newcommand{\Excom}{E_{\ranx\co}}

\newcommand{\varparax}{\var_{\para,\ranx}}

\newcommand{\bern}{\mbox{Bernoulli}}
\newcommand{\nuisprior}{\psi}
\newcommand{\nuis}{\beta}
\newcommand{\nuisnull}{\nuis_0}

\newcommand{\nuisnulllin}{\nuis_{0,\text{lin}}}
\newcommand{\nuisnullquad}{\nuis_{0,\text{quad}}}
\newcommand{\nuisnullcubic}{\nuis_{0,\text{cubic}}}
\newcommand{\nuispriornull}{\psi_0}
\newcommand{\nsupnull}{B_0}
\newcommand{\nuisalt}{\beta_1}
\newcommand{\nuisaltinthatj}{\hat{\nuis}_{1,\text{int}}^{(j)}}
\newcommand{\nuisaltslopehatj}{\hat{\nuis}_{1,\text{slope}}^{(j)}}
\newcommand{\nuisnullinthatj}{\hat{\nuis}_{0,\text{int}}^{(j)}}
\newcommand{\nuisnullslopehatj}{\hat{\nuis}_{0,\text{slope}}^{(j)}}
\newcommand{\nuisinthatji}{\hat{\nuis}_{i,\text{int}}^{(j)}}
\newcommand{\nuisslopehatji}{\hat{\nuis}_{i,\text{slope}}^{(j)}}
\newcommand{\nuisprioralt}{\psi_1}

\newcommand{\altmeanint}{\altmean_{\text{int}}}
\newcommand{\altmeanslope}{\altmean_{\text{slope}}}
\newcommand{\nuisnullint}{\nuis_{0,\text{int}}}
\newcommand{\nuisnullslope}{\nuis_{0,\text{slope}}}
\newcommand{\nuisaltint}{\nuis_{1,\text{int}}}
\newcommand{\nuisaltslope}{\nuis_{1,\text{slope}}}

\newcommand{\nuisinttheta}{\nuis_{\para,\text{int}}}
\newcommand{\nuisslopetheta}{\nuis_{\para,\text{slope}}}
\newcommand{\nuisinti}{\nuis_{i,\text{int}}}
\newcommand{\nuisslopei}{\nuis_{i,\text{slope}}}

\newcommand{\nsupalt}{B_1}

\newcommand{\nsup}{B}
\newcommand{\nsalt}{B_1}

\newcommand{\postv}{posterior-prior ratio }
\newcommand{\logv}{log }

\newcommand{\postbayes}{P_{\text{bayes}}}
\newcommand{\postlrt}{P_{\text{LRT}}}
\newcommand{\tkbayes}{TK_{\text{bayes}}}
\newcommand{\tklrt}{TK_{\text{LRT}}}

\newcommand{\mpower}{M_{\text{power}}}
\newcommand{\mspread}{M_\text{spread}}

\tikzset{
    state/.style={
           rectangle,
           rounded corners,
           draw=black, very thick,
           minimum height=2em,
           inner sep=4pt,
           text centered,
           execute at begin node=\setlength{\baselineskip}{1.5em}
           },
                          descr/.style={
        fill=white,
        inner sep=2.5pt
    },
    connector/.style={
     -latex,
     font=\scriptsize
    },
    rectangle connector/.style={
        connector,
        to path={(\tikztostart) -- ++(#1,0pt) \tikztonodes |- (\tikztotarget) },
        pos=0.5,
        rounded corners
    },
    rectangle connector/.default=-2cm,
    straight connector/.style={
        connector,
        to path=--(\tikztotarget) \tikztonodes
    }
}

\definecolor{darkblue}{rgb}{0.0,0.0,0.5}
\RequirePackage[colorlinks,allcolors=darkblue]{hyperref}

\newcommand{\blind}{1}

\begin{document}

\def\spacingset#1{\renewcommand{\baselinestretch}%
{#1}\small\normalsize} \spacingset{1}


\if0\blind
{

} \fi

\if1\blind
{
  \title{\bf Designing Test Information\\and Test Information in Design}
  \author{David E. Jones$^1$ and Xiao-Li Meng$^2$\\\\\hspace{.2cm}Texas A\&M University$^1$ and Harvard University$^2$\hspace{.2cm}\\\\
    }
  \maketitle
} \fi
\kern -.2in
\begin{abstract}
\citet{degroot1962} developed a general framework for constructing Bayesian measures of the expected information that an experiment will provide for estimation.  We propose an analogous framework for measures of information for hypothesis testing. In contrast to estimation information measures that are typically used  for surface estimation, test information measures are more useful in experimental design for hypothesis testing and model selection. In particular, we obtain a probability based measure, which has more appealing properties than variance based measures in design contexts where decision problems are of interest. The underlying intuition of our design proposals is straightforward: to distinguish between models we should collect data from regions of the covariate space for which the models differ most.   \citet{frac_info} gave an asymptotic equivalence between their test information measures and Fisher information. We extend this result to all test information measures under our framework. Simulation studies and an application in astronomy demonstrate the utility of our approach, and provide comparison to other methods including that of \citet{boxhill1967}.
\end{abstract}

\noindent%
{\it Keywords:}  Statistical information; Optimal design; Bayes factors; Hypothesis testing; Model selection; Power

\spacingset{1.45} 

\section{Motivation and Overview}
\label{sec:overview}

\subsection{Test information framework: foundations and developments}
\label{sec:intro}

\cite{frac_info} highlighted that the amount of information provided by an experiment depends on our goals, and in particular the amount of information for hypothesis testing can be different from that for estimation. Nonetheless, the importance of information measures and the need for a framework for constructing and understanding them is common to both the testing and estimation scenarios. Indeed, \citet{ginebra2007} emphasized that
flexible information measures are essential because information is a ``highly
multi-dimensional concept" that cannot be adequately captured by a narrow definition.
In statistics, the most common appeal to information is the use of Fisher information to determine the asymptotic variance of the maximum likelihood estimator (MLE). However, the key to the importance of information measures is that they quantify what it is {\it possible} to learn on average given a data generating model (for the data to be used), and thus they go beyond detailing the properties of a specific procedure. This is illustrated by the fact that the Fisher information is not {\it merely} related to the asymptotic variance of the MLE, it also appears in the Cram\'{e}r-Rao lower bound for the variance of all unbiased estimators. Given such requirements, it is natural that for Bayesian inference, estimation information measures should be based on the posterior distribution, which captures  all available information. Furthermore, it should be no surprise that the fundamental component of our test information measures is the Bayes factor or likelihood ratio.

General development of the existing estimation information framework began with the concept of statistical entropy introduced by \cite{shannon}:
\begin{align}
H(\prior) = \Epara[-\log\prior (\para)] = -\int_\psup \prior(\para)\log \prior(\para) \mu(d\para).
\label{eqn:entropy}
\end{align}
Entropy measures the information gained by observing the random variable $\para \in \psup $ with density $\prior$, with respect to the measure $\mu$. \cite{lindley1956} defined the expected information about a parameter provided by an experiment as the difference between the prior entropy and the expected posterior entropy. This measure has seen many applications, such as defining D-optimality  \citep[see the review by][]{chaloner} and reference priors \citep[e.g.,][]{berger2009}. \cite{degroot1962} generalizes Lindley's framework by replacing Shannon's entropy by a general measure of {\it uncertainty}.

Our first contribution is to synthesize the test information measures suggested by \cite{frac_info} to develop an analogous framework to that of \cite{degroot1962} for test information. The general measures of expected test information that we propose use the $f$-divergence introduced by \cite{csiszar1963} and \cite{ali_silvey1966}, and we extend the concept to define observed and conditional test information because these are of great importance in sequential design. This construction and the parallels with estimation information are summarized in Figure \ref{fig:theory_diagram}. The testing and estimation scenarios differ because the hypotheses to be compared must be incorporated into test information measures, and in computing expected test information we must choose a hypothesis to condition on. This last feature suggests that every expected test information measure should have a dual which conditions on the other hypothesis. Previous authors, including \citet{lopez}, have identified a similar phenomenon, but have not resolved the issue of two differing measures. Here, we demonstrate that there exists an appealing subclass of our measures which give the same expected test information as their duals, thereby resolving the conflict.

Our second contribution is to establish further connections between test and estimation information measures. These connections concern an important quantity for sequential design discussed by \cite{nicolae2004} and \cite{frac_info}, namely the fraction of information contained in the observed data relative to that contained in the intended complete data, which is not fully observed. \cite{frac_info} established an asymptotic equivalence between their measures of the fraction of observed test information and the fraction of observed Fisher information (for estimation), as the distance between the null parameter and the MLE goes to zero. We show that, by allowing different weighting of observed and missing Fisher information, the equivalence can be extended to hold for all test information measures under our framework. This result
identifies an appealing class of test information measures that weight observed and missing {\it estimation} information equally (in the limit considered).

With the basic foundations of our test information framework in place, we consider its practical implications. \cite{nicolae2004} and \cite{frac_info} put forward their  measures of the fraction of observed (or missing) test information with the purpose of informing data collection decisions in genetic linkage studies. We now build on this by illustrating specifically how test information measures may be used in experimental design, both in model selection and coefficient testing scenarios. In the design for model selection scenario, it is often not clear how to use estimation information measures, but the use of test information measures is intuitive. We demonstrate this advantage by finding optimal designs for choosing between the complementary log-log and Probit link functions for binary linear regression.

Next, in the specific case of testing for linear regression parameters, we give a closed form design criterion that is related to the familiar Bayesian alphabet optimality criteria, and demonstrate its use in sequential design. We also propose a posterior probability based expected test information measure, which has many appealing properties, including the aforementioned duality (and can be applied more widely than the linear regression context). Our approach to design for testing linear regression parameters differs from that of \citet{toman1996}, which minimizes Bayes risk, because we instead maximize the expected probabilistic information for distinguishing hypotheses.

Among the existing literature on measuring design information
for hypothesis testing and model selection, the approach of \citet{boxhill1967} is perhaps most similar to ours.
They chose designs that maximize the expected change in entropy of the posterior probability mass function on the model indicator, but did not provide a general framework for test information measures and design.  
In terms of mathematical justification, our framework benefits from the work of \citet{ginebra2007}, which reviewed and synthesized previous theory including that of \citet{blackwell1951}, \citet{sherman1951}, \citet{stein1951}, and \citet{lecam1964}. Specifically, our expected test information measures satisfy (up to aesthetics) the three basic requirements identified in \citet{ginebra2007}.
Our framework adds the concept of coherent dual test information measures, and our observed test information measures have fewer restrictions than those suggested by \citet{ginebra2007} (who focused on the estimation case). 
Furthermore, we demonstrate that, in practice, test and estimation information measures behave very differently, despite their common mathematical roots described by \citet{ginebra2007}.

Other important design approaches for  model discrimination include the $D_s$-optimal designs introduced by \citet{stigler1971optimal} and \citet{studden1980d_s}, and the $T$-optimal designs introduced by \citet{atkinson1975design} and \citet{atkinson1975optimal}. The former approach focuses on nested models and the latter on the squared differences between mean functions, and both differ somewhat to our more general framework. The KL-optimality criterion of \citet{lopez} is closely related to the likelihood ratio test and therefore some specific measures falling under our framework, though we principally consider the Bayes factor.  Connections between T-optimal, $D_s$-optimal, and KL-optimal designs have been established by a number of authors including \citet{dette2009optimal} and \citet{lopez}. However, connections between designs for estimation and these model discrimination designs are rarer, partly because optimal designs for the distinct goals of estimation and model discrimination are usually different.

A key limitation of designs optimized for distinguishing between a set of models is their inherent sensitivity to these hypothesized models. In particular, if none of the hypothesized models reasonably capture reality, then the design selected by any model-based measure can be seriously sub-optimal. Furthermore, designs for model discrimination tend to have only a small number of unique design points and therefore model checking beyond the intended discrimination  can be difficult, e.g., it may not even be possible to estimate the parameters of the hypothesized models under the optimal design, see for instance \citet{dette2009optimal}. These issues are almost unavoidable and can perhaps only be mitigated by generic space filling designs sometimes used in estimation scenarios for similar reasons, with the familiar trade-off between robustness and efficiency. With that in mind, test information measures are valuable design tools when scientific investigations are of a confirmatory nature, meaning that there is some reason to believe that the proposed models adequately capture accepted scientific phenomena (physics is one area where this quite often holds). In summary, the information measures we propose are beneficial whenever an investigator seeks to compare several reasonable competing models.

Our paper is organized as follows. The remainder of Section \ref{sec:overview} gives three categories of scientific problems where test information measures are useful, briefly reviews the estimation information framework proposed by \cite{degroot1962}, and discusses the parallels with the test information measures introduced by \cite{frac_info}. The main body of the paper is divided between Sections \ref{sec:expected} and \ref{sec:observed}, which deal with expected and observed test information, respectively. These sections finish with illustrations of the practical use of test information in design and sequential design for decision problems, respectively. Section \ref{sec:expected} also introduces a fundamental symmetry condition that defines an appealing class of test information measures. Section \ref{sec:theory} presents our main result linking test and estimation information. In Section \ref{sec:application} we apply our approach to a problem in astronomy, and Section  \ref{sec:discussion} concludes with discussion and open challenges.

\subsection{Uses of test information}
\label{sec:uses}

We now describe three experimental design problems, representing broad categories of scientific questions.

{\it Classification and model selection.} In astronomy, the intensity of some sources (e.g., Cepheid stars) varies periodically over time, thus creating a continuous function of time called a lightcurve. Such sources can be classified by features of their lightcurves, e.g., the period or characteristic shape. Since telescope time is limited for any group of researchers, the lightcurve of a source is observed at a number of time points and then a classifier is applied. For example, some modern techniques use random forest classifiers, e.g., \cite{richards2011} and \cite{dubath2011}. 
Intuitively, given some lightcurve observations, the design problem is to pick the time of the next observation that will maximize the probability of correct classification.

{\it Screening and follow-up.} In genetic linkage studies it is of interest to test if markers (or genes) located close together on the same chromosome are more likely to be inherited together than if markers are inherited independently (the null hypothesis). This is a screening process because the magnitude of the linkage (i.e., dependence) is ultimately of interest. 
In the case of too much missing information, a follow-up study could choose between increasing the number of markers in potential regions of linkage or collecting samples from additional individuals. To assess which option is likely to have greater power, for example, we must take the models under the two hypotheses into account, as studied in \cite{nicolae2004} and \cite{frac_info}.

{\it Robust design.} Test information measures can also be useful in applications at the interface of testing and estimation. In chemical engineering, it is often of interest to estimate the mean yield of a product under different conditions, and ultimately to model the yield. In this situation, space filling designs are usually preferred because it is unknown where the regions of high (and low) yield will be. However, space filling designs can vary in their efficiency for distinguishing specific models, and test information measures  can be used to select the ones that best separate important candidate models.

\subsection{Bayesian information for estimation}
\label{sec:degroot}

We briefly review the framework of \cite{degroot1962} to help make clear both distinctions and parallels between test and estimation information. Suppose that we are interested in a parameter $\para \in \psup$ and our prior distribution is $\prior$. Information about $\theta$ is gained through an experiment $\ept$ whose future outcome is $\ranx \in \xsup$, the set of possible outcomes of the experiment. For example, $\ept$ may specify the design points at which data are to be collected. We denote by $\info(\ept;\prior)$ a measure of the expected information to be gained by conducting $\ept$.
The measure $\info(\ept;\prior)$ should have basic properties such as non-negativity and additivity. To specify the meaning of additivity we need the notion of conditional information: if $\ept=(\ept_1,\ept_2)$ is an experiment composed of two sub-experiments, then we denote by $\info(\ept_2|\ept_1;\prior)$ the expected {\it conditional} information to be gained by conducting $\ept_2$ {\it after} conducting $\ept_1$, i.e., the expected new information that will be gained from $\ept_2$. Additivity can now be formalized.

\begin{defi}
\label{def:add}
An information measure $\info$ is additive if, for any composite experiment $\ept=(\ept_1,\ept_2)$ and any proper prior $\prior$, the following relation holds
\begin{eqnarray}
\label{eqn:add}
\info(\ept;\prior) = \info(\ept_1;\prior)+\info(\ept_2|\ept_1;\prior).
\end{eqnarray}
\end{defi}

\cite{degroot1962} chose $\info(\ept;\prior)$ to be the difference between the prior uncertainty and the expected posterior uncertainty about $\para$. He defined the prior uncertainty to be $\uncp(\prior)$, where the {\it uncertainty function} $\uncp$ is a concave functional of $\prior$, i.e., $\uncp(\lambda \prior_1 +(1-\lambda) \prior_2) \geq \lambda\uncp( \prior_1) +(1-\lambda)\uncp( \prior_2)$ for any two densities $\prior_1$ and $\prior_2$ and $\lambda\in [0,1]$. Similarly, \cite{degroot1962}  defined the expected posterior uncertainty to be $\Ex[\uncp(p(\cdot|\ranx))]$, where the expectation is with respect to $f(\datx) = \int_\psup f(\datx|\para)\prior(\para)\mu(d\para)$. Thus, we have the following.
\begin{defi}
\label{def:est_info}
The expected Bayesian estimation information provided by an experiment $\ept$, under a proper prior $\prior$, is
\begin{eqnarray}
\label{eqn:beinfo}
\info(\ept;\prior) = \uncp(\prior) - \Ex[\uncp(p(\cdot|\ranx))].
\end{eqnarray}
\end{defi}

\cite{lindley1956} suggested $\uncp$ should be the entropy function $H$ given in (\ref{eqn:entropy}). 
\cite{degroot1962} showed that (\ref{eqn:beinfo}) satisfies non-negativity for all priors $\prior$ and experiments $\ept$ if and only if $\uncp$ is concave. To generalize further, we follow 
the logic of Definition \ref{def:est_info} and define the expected conditional estimation information contained in the second of two sub-experiments as
\begin{eqnarray}
\label{eqn:con_beinfo}
\info(\ept_2|\ept_1;\prior) =  \Exone[\uncp(p(\cdot|\ranx_1))] - \Ex[\uncp(\post(\cdot|\ranx))],
\end{eqnarray}
where $\ranx=\{\ranx_1, \ranx_2\}$,  and $\ranx_i$ is
the outcome of $\ept_i$, $i=1,2$. The desired additivity then follows trivially by definition. However, it is not generally true that $\info(\ept;\prior)=\info(\ept_1;\prior)+\info(\ept_2;\prior)$, even when $\ranx_1$ and $\ranx_2$ are independent.

\subsection{Test information measures proposed by \citet{frac_info}}
\label{sec:parallels}

\def\twidth{5.75cm}
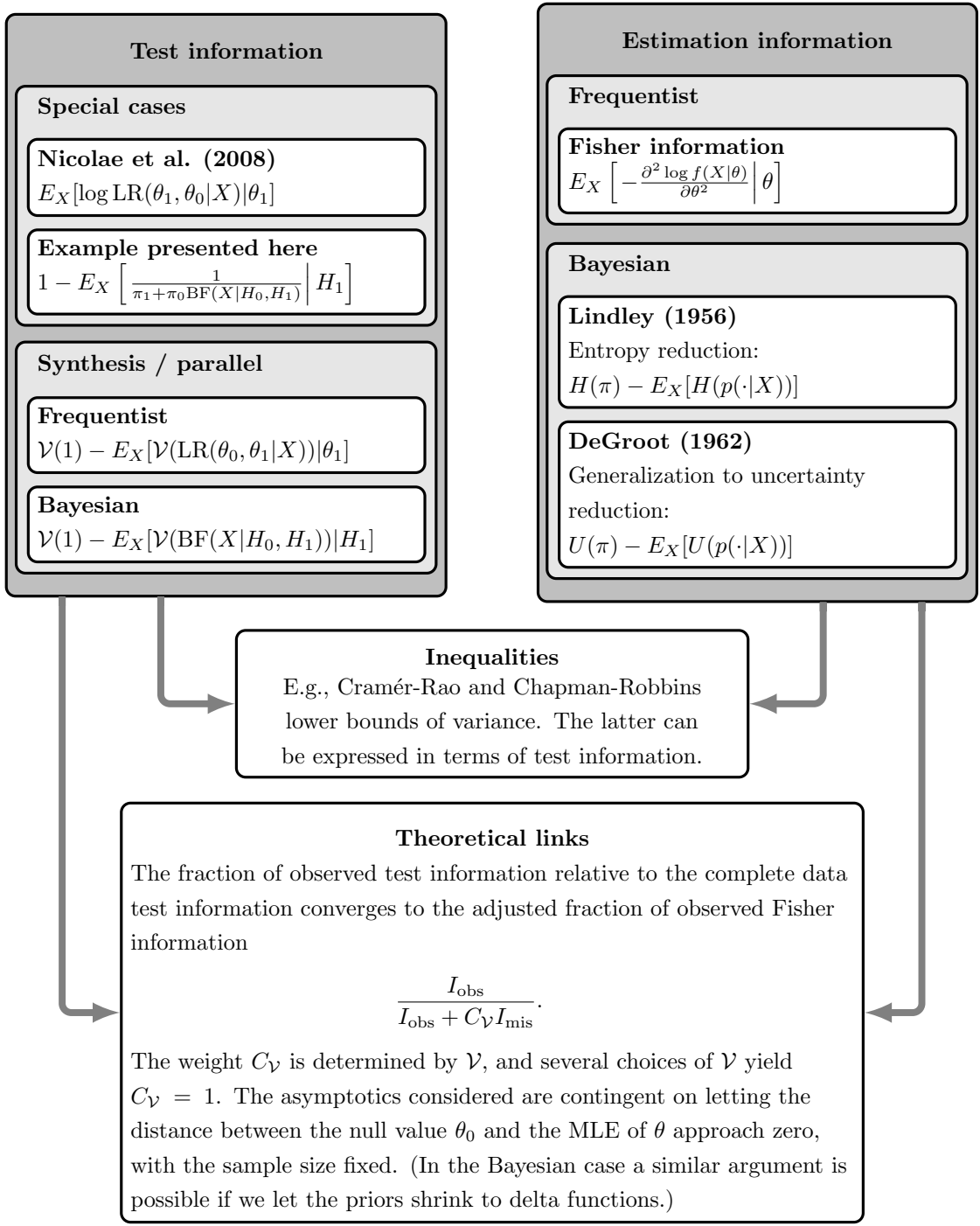
\begin{figure}[p]
\begin{center}
\begin{tikzpicture}


\node[state, fill=white, text width=\twidth, align =left] (nicolae){
    \textbf{Nicolae et al. (2008)}\\
    $\Ex[\log \text{{\usefont{OT1}{cmr}{l}{sc} LR}}(\para_1,\para_0|\ranx)|\para_1]$};

\node[state, fill=white, text width=\twidth,align =left,below = 0.5em of nicolae] (postef){
    \textbf{Example presented here}\\
    $1 - \Ex\left[\left.\frac{1}{\pi_1+\pi_0 \text{{\usefont{OT1}{cmr}{l}{sc} BF}}(\ranx|H_0,H_1)}\right|H_1\right]$};

\node[state, fill=white, text width=\twidth,align =left,below = 3.5em of postef] (testfreq){
    \textbf{Frequentist}\\
    $\decz - \Ex[\declr(\text{{\usefont{OT1}{cmr}{l}{sc} LR}}(\para_0,\para_1|\ranx))|\para_1]$};

\node[state, fill=white, text width=\twidth,align =left,below = 0.5em of testfreq] (testbayes){
    \textbf{Bayesian}\\
    $\decz - \Ex[\declr(\text{{\usefont{OT1}{cmr}{l}{sc} BF}}(\ranx|H_0,H_1))|H_1]$};

\node[above = 0.5em of nicolae,align =left,text width = \twidth] (testcase){
\textbf{Special cases}};

\node[above = 1em of testcase] (test) {
    \textbf{Test information}};

\node[above = 0.5em of testfreq,align =left,text width = \twidth] (synthesis){\textbf{Synthesis / parallel}};


\node[state, fill=white, text width=\twidth, align =left,right = 2cm of nicolae,yshift=0.1cm] (fisher){
    \textbf{Fisher information}\\
    $\Ex\left[\left.-\frac{\partial^2\log f(\ranx|\para)}{\partial \para^2}\right|\para\right]$};

\node[state, fill=white, text width=\twidth,align =left,below = 3.5em of fisher] (lindley){
    \textbf{Lindley (1956)}\\
    Entropy reduction:\\
    $H(\prior) - \Ex[H(p(\cdot|\ranx))]$};

\node[state, fill=white, text width=\twidth,align =left,below = 0.5em of lindley] (degroot)
    {\textbf{DeGroot (1962)}\\
    Generalization to uncertainty\\ reduction:\\
    $\uncp(\prior) - \Ex[\uncp(p(\cdot|\ranx))]$};

\node[above = 0.5em of fisher,align =left,text width = \twidth] (estfreq){
\textbf{Frequentist}};

\node[above = 1em of estfreq] (estimation) {
    \textbf{Estimation information}};

\node[above = 0.5em of lindley,align =left,text width = \twidth] (estbayes){\textbf{Bayesian}};


    \begin{pgfonlayer}{background}
        \path (fisher.west |- estimation.north)+(-0.3,0.3) node (a) {};
        \path (degroot.south east)+(+0.3,-0.5) node (b) {};
        \path[fill=gray!50!white,rounded corners, draw=black, very thick] (a) rectangle (b);

        \path (estbayes.north west)+(-0.2,0.05) node (a) {};
        \path (degroot.south east)+(+0.15,-0.15) node (b) {};
        \path[fill=gray!20!white,rounded corners,  draw=black, very thick] (a) rectangle (b);

        \path (estfreq.north west)+(-0.2,0.05) node (a) {};
        \path (fisher.south east)+(+0.15,-0.15) node (b) {};
        \path[fill=gray!20!white,rounded corners, draw=black, very thick] (a) rectangle (b);
    \end{pgfonlayer}

    \begin{pgfonlayer}{background}
        \path (testcase.west |- test.north)+(-0.35,0.3) node (a) {};
        \path (testbayes.south east)+(+0.3,-0.5) node (b) {};
        \path[fill=gray!50!white,rounded corners, draw=black, very thick] (a) rectangle (b);

        \path (testcase.north west)+(-0.2,0.05) node (a) {};
        \path (postef.south east)+(+0.15,-0.15) node (b) {};
        \path[fill=gray!20!white,rounded corners,  draw=black, very thick] (a) rectangle (b);

        \path (synthesis.north west)+(-0.2,0.05) node (a) {};
        \path (testbayes.south east)+(+0.15,-0.15) node (b) {};
        \path[fill=gray!20!white,rounded corners, draw=black, very thick] (a) rectangle (b);
    \end{pgfonlayer}

\node[state,below=1cm of testbayes,xshift=4.04cm,text width=7.5cm] (bounds){
    \begin{center}\vspace{-0.15cm}\textbf{Inequalities}\end{center}
    \vspace{-0.075cm}
    E.g., Cram\'{e}r-Rao and Chapman-Robbins lower bounds of variance. The latter can be expressed in terms of test information.};

\node[state,below=1em of bounds,align=left,text width=11cm] (reconciliation){
    \begin{center}\vspace{-0.3cm}\textbf{Theoretical links}\end{center}
    \vspace{-0.275cm}
    The fraction of observed test information relative to the complete data test information converges to the adjusted fraction of observed Fisher information
    \begin{eqnarray*}
    \frac{I\ob}{ I\ob + \cn I\mis}.
    \end{eqnarray*}
    The weight $\cn$ is determined by $\declr$, and several choices of $\declr$ yield $\cn=1$. The asymptotics considered are contingent on letting the distance between the null value $\para_0$ and the MLE of $\para$ approach zero, with the sample size fixed. (In the Bayesian case a similar argument is possible if we let the priors shrink to delta functions.)};

\draw [rectangle connector=0cm,line width=3pt,gray] (degroot.south)+(1,-0.65) to (bounds.east);
\draw [rectangle connector=0cm,line width=3pt,gray] (degroot.south)+(2.5,-0.65) to (reconciliation.east);
\draw [rectangle connector=0cm,line width=3pt,gray] (testbayes.south)+(-1,-0.65) to (bounds.west);
\draw [rectangle connector=0cm,line width=3pt,gray] (testbayes.south)+(-2.5,-0.65) to (reconciliation.west);
\end{tikzpicture}
\end{center}
\caption{Summary of estimation and test information theory. The synthesis of test information measures into one coherent framework paralleling the estimation framework is new. Also new are the general links between estimation and test information, although \cite{frac_info} considered the same connection with Fisher information for specific cases. \label{fig:theory_diagram}}
\end{figure}
For the sharp test hypotheses $H_0:\para=\para_0$ and $H_1:\para=\para_1$, \cite{frac_info} (implicitly) proposed the  very natural frequentist expected test information measure
\begin{eqnarray}
\info^T(\ept;\para_0,\para_1) = \Ex[\log \lr(\para_1,\para_0|\ranx)|\para_1],
\label{eqn:log_measure}
\end{eqnarray}
where the superscript $T$ indicates the testing context, and
\begin{eqnarray}
\label{eqn:lr}
\lr(\para_0,\para_1|x) = \frac{f(\datx|\para_0)}{f(\datx|\para_1)}
\end{eqnarray}
is the likelihood ratio. We observe that (\ref{eqn:log_measure}) is the Kullback-Leibler (KL) divergence between the data models $f(\cdot|\para_0)$ and $f(\cdot|\para_1)$, and thus it is closely connected to the entropy based measure proposed by \cite{lindley1956}. (The KL divergence between two densities $g$ and $h$, will be denoted $KL(g||h)$, and is defined as $\int_\xsup g(x)\log\left(g(x)/h(x)\right) \mu(dx)$, where the support of $g$ is assumed to be a subset of the support of $h$.) Nonetheless, there is a good reason why \cite{frac_info} did not simply use Definition \ref{def:est_info} to construct measures of test information; namely, it does not take the hypotheses into account. Indeed, the presence of the two parameter values, $\para_0$ and $\para_1$, in (\ref{eqn:log_measure}) clearly distinguishes test information from the estimation information we have considered so far. This difference makes intuitive sense because it represents the difference between gaining evidence for distinguishing two hypotheses and neutrally gaining knowledge about the parameter.

In practice, the alternative hypothesis is often composite and in the Bayesian context we then write $H_1: \para \sim \prior$, for some prior $\prior$. One of the Bayesian measures of expected test information (implicitly) suggested by \cite{frac_info} is
\begin{eqnarray}
\info^T(\ept;\para_0,\prior) = \varparax(\lr(\para_0,\para|\ranx)).
\label{eqn:var_measure}
\end{eqnarray}
Variance and entropy are both measures of spread and hence (\ref{eqn:var_measure}) is also connected to the measure proposed by \cite{lindley1956}, although no logarithm is taken in (\ref{eqn:var_measure}). The key distinction with Definition \ref{def:est_info} is again due to the appearance of the null hypothesis $\para_0$.  In summary, these examples have connections with the estimation information measures reviewed in Section \ref{sec:degroot}, but also have common features distinguishing test information from estimation information. Based on these parallels and distinctions, the next section proposes our general framework for constructing test information measures, which is summarized in the grey box on the left of Figure \ref{fig:theory_diagram}.

\section{Expected Test Information: Theory and Applications}
\label{sec:expected}

\subsection{Test information: a synthesis}
\label{sec:expected_test}

The two key properties of expected information measures are non-negativity and additivity. For simplicity, we develop our framework in the case of continuous densities and the Lebesgue measure. Theorem 2.1 of \cite{degroot1962} establishes non-negativeness of the estimation information reviewed in Section \ref{sec:degroot}.  Writing the marginal density of $\datx$ as $f(\datx) = \int_\psup f(x|\para)\prior(\para) d\para$, the key formula underlying the theorem is
\begin{eqnarray}
\Ex[p(\cdot|\ranx)] = \int_\xsup p(\cdot|\datx)f(\datx)d\datx = \prior(\cdot).\label{eqn:degroot_jensen}
\end{eqnarray}
That is, the expected posterior density with respect to the marginal density is the prior density. To see the corresponding key identity for hypothesis testing, we first observe that the expected test information (\ref{eqn:log_measure}) uses the likelihood ratio as the fundamental statistic for quantifying the information for distinguishing two values of $\para$. More generally, the hypotheses may be composite, say $H_0: \para \in \psupnull$ and $H_1: \para \in \psupalt$, in which case we turn to the Bayesian perspective and replace the likelihood ratio with the Bayes factor
\begin{align}
\label{eqn:bf}
\bfac(\datx|H_0,H_1) = \frac{f(\datx|H_0)}{f(\datx|H_1)} = \frac{\int_\psupnull f(\datx|\para)\prior(\para|H_0)d\para}{\int_\psupalt f(\datx|\para)\prior(\para|H_1) d\para} = \frac{\int_\psupnull f(\datx|\para)\frac{\prior(\para)}{\prior_0}d\para}{\int_\psupalt f(\datx|\para)\frac{\prior(\para)}{\prior_1}d\para},
\end{align}
where $\prior_i=P(\para \in \psupi)$, for $i=0,1$. (We assume $\prior_i \not=0$ throughout, for $i=0,1$.) Thus, for hypothesis testing, the analogous formula to (\ref{eqn:degroot_jensen}) is
\begin{eqnarray}
\label{eqn:like1}
\Ex[\bfac(\ranx|H_0,H_1)|H_1] = \int_\xsup \frac{f(\datx|H_0)}{f(\datx|H_1)}f(\datx|H_1)d\datx = 1.
\end{eqnarray}
That is, the expected Bayes factor (or likelihood ratio), under the alternative, does not favor either hypothesis. For simplicity, we assume here and throughout that the support of $f(\cdot|\para)$ is $\xsup$ for all $\para \in \psup$. Equation (\ref{eqn:like1}) allows us to apply Jensen's inequality to ensure that the general expected test information given in Definition \ref{def:decinfo} (below) is non-negative. For test information, the parallel of the uncertainty function $\uncp$ is the {\it evidence function} $\declr$, which acts on the positive real numbers and in particular has the Bayes factor (or likelihood ratio) as its argument. The use of Jensen's inequality to ensure non-negativity requires that the evidence function is concave, and we therefore assume concavity throughout. Note that what is measured by the evidence function is the evidence {\it in support of the null hypothesis}, and therefore, like \cite{degroot1962}, we are interested in a reduction, i.e., a reduction in the evidence supporting the null hypothesis. 

\begin{defi} \label{def:decinfo} Under $H_1:\para \in \psupalt$, the expected test information provided by the experiment $\ept$ for comparing the hypotheses $H_0: \para \in \psupnull$ and $H_1: \para \in \psupalt$, for a given evidence function $\declr$ and a proper prior $\prior$, is defined as
\begin{align}
\bdinfofull(\ept;H_0,H_1,\prior) &= \decz - \Ex[\declr(\bfac(\ranx|H_0,H_1))|H_1],
\label{eqn:bayes_decinfo}
\end{align}
where $\psupnull \cap \psupalt = \emptyset$.
\end{defi}

The prefix ``expected" is necessary because the Bayesian approach generally assumes that data have been observed. Note that (\ref{eqn:bayes_decinfo}) is mathematically equivalent to the frequentist measure
\begin{align}
\decz - \Ex[\declr(\lr(\para_0,\para_1|\ranx)|\para_1] \label{eqn:freq_decinfo}
\end{align}
when $\psupi=\{\para_i\}$ and $P(\para = \para_i)\not=0$, for $i=0,1$. (The frequentist perspective is also recovered if the prior is viewed as part of the data generating model.) Under (\ref{eqn:freq_decinfo}), the measure (\ref{eqn:log_measure}) is given by choosing $\declr(\zlr)=\log(\zlr)$. The mathematical equivalence of Bayesian and frequentist measures of expected test information means that we can interchange the Bayes factor in (\ref{eqn:bayes_decinfo}) and the likelihood ratio as convenient. More generally, the Bayes factor in (\ref{eqn:bayes_decinfo}) can be replaced by {\it any} numerical comparison of the hypotheses, at least if the baseline is also adjusted. However, the main focus here will be on the Bayesian perspective because it is statistically coherent and is conceptually well suited to incorporating composite hypotheses (and nuisance parameters, see Section \ref{sec:nuisance}) when no data have been observed, as is often the case when we choose a design. We retain the argument $\prior$ in our notation $\bdinfofull(\ept;H_0,H_1,\prior)$ as a reminder that (\ref{eqn:bayes_decinfo}) does depend on the prior $\prior$, which we should therefore choose carefully, as with the specification of any part of our models. Also note that the parameter in Definition \ref{def:decinfo} can simply be a model indicator and hence our framework goes beyond parametric models. 

The final term of (\ref{eqn:bayes_decinfo}) is the $f$-divergence introduced by \cite{csiszar1963} and \cite{ali_silvey1966}, which generalizes KL divergence. Indeed, as mentioned in Section \ref{sec:parallels}, the measure  (\ref{eqn:log_measure}) is a KL divergence. The properties of KL divergence alert us to the important feature that expected test information is not necessarily symmetric in the two hypotheses. A class of evidence function that treat the hypotheses equally will be introduced in Section \ref{sec:symmetry_class}. 


The baseline term $\declr(1)$ ensures non-negativity of expected test information, and has common intuitive appeal because it represents no evidence for either hypothesis. Note that we have chosen a baseline that does not incorporate investigator specific quantities such as the frequentist size of a particular test, because a general information measure should appeal to many different researchers. Information measures with general appeal are typically {\it maximal information} measures, for example, Fisher information measures the maximal estimation information asymptotically available. Our test information measures are also implicitly maximal since they are functions of only the prior probabilities and Bayes factor (or likelihood ratio), which contain all the relevant information. For example, the KL divergence in (\ref{eqn:log_measure}) uses the expected log likelihood ratio to quantify the maximal probabilistic information for distinguishing the hypotheses, and leaves the choice of a specific decision rule to individual investigators. 

From this point on, we will frequently write $\bdinfo(\ept)$ to mean $\bdinfofull(\ept;H_0,H_1,\prior)$ and for similar notation will again drop the arguments after the semicolon when this causes no confusion. We complete our initial development of test information by specifying the form of $\bdinfo(\ept_2|\ept_1;H_0,H_1,\prior)$, which is easily deduced from the expected conditional estimation information (\ref{eqn:con_beinfo}).

\begin{defi}
\label{def:con_decinfo}
The expected conditional test information $\bdinfo(\ept_2|\ept_1)$ provided by conducting the experiment $\ept_2$ after $\ept_1$ is
\begin{align}
\label{eqn:con_decinfo}
\bdinfofull(\ept_2|\ept_1;H_0,H_1,\prior)= \Exone[\dec(\ranx_1)|H_1] - \Ex[\dec(\ranx)|H_1],
\end{align}
where $\dec(Z) = \declr(\bfac(Z|H_0,H_1))$, for $Z=\ranx_1$
and $Z=\ranx=\{\ranx_1,\ranx_2\}$.
\end{defi}

That (\ref{eqn:con_decinfo}) is non-negative is again a consequence of Jensen's inequality:
\begin{align}
\Extwo[\dec(\datx_1,\ranx_2)|H_1,\datx_1] \leq  \dec(\datx_1),
\label{eqn:con_jensen}
\end{align}
where, to make the expressions easier to read, we have denoted the observed data by lower case letters, and unobserved data by upper case letters. Given Definition \ref{def:con_decinfo}, the additivity property of Definition \ref{def:add} holds trivially, i.e., $\bdinfo(\ept) = \bdinfo(\ept_1)+\bdinfo(\ept_2|\ept_1)$.

\subsection{Symmetry class and a probability based measure}
\label{sec:symmetry_class}

The best choice of $\declr$ will to some extent depend on the particular context (see Section \ref{sec:theory} for some theoretical guidance), but here we propose a class of evidence functions that have appealing properties. The class is those evidence functions that treat the hypotheses symmetrically and in particular satisfy the condition
\begin{align}
\frac{\declr(\zlr;H_0,H_1)}{\declr(\zlr^{-1};H_1,H_0)}=\zlr.
\label{eqn:observe_coherence}
\end{align}
Naturally, $\declr(\zlr^{-1};H_1,H_0)$ represents the evidence for the alternative, since the roles of $H_0$ and $H_1$ have been swapped. Thus, setting $z = \bfac(\datx|H_0,H_1)$, the symmetry condition (\ref{eqn:observe_coherence}) states that our choice of $\declr$ should preserve the Bayes factor. We include the arguments $H_0$ and $H_1$ in (\ref{eqn:observe_coherence}) because in general the evidence measures may be allowed to depend on the order of the hypotheses through prior probabilities as well as through the Bayes factor.

The symmetry condition (\ref{eqn:observe_coherence}) implies that the resulting expected test information measure satisfies the fundamental coherence identity
\begin{align}
\bdecinfo(\ept;H_0,H_1,\prior)=\bdecinfo(\ept;H_1,H_0,\prior).
\label{eqn:info_link}
 \end{align}
The right hand of (\ref{eqn:info_link}) swaps the hypotheses, indicating that the {\it dual} expected test information measure $\bdecinfo(\ept;H_1,H_0,\prior)$ takes an expectation with respect to $f(\cdot|H_0)$, rather than $f(\cdot|H_1)$. Indeed, the dual test information measure quantifies the reduction in evidence for the alternative when data are collected under the null. The coherence identity (\ref{eqn:info_link}) states that, before we observe any data, the expected amount of information in the data for distinguishing the two hypotheses is the same regardless of which is in fact true. This symmetry is intuitive because the probabilistic separation of the two marginal data models $f(\cdot|H_0)$ and $f(\cdot|H_1)$ should not depend on which hypothesis is true. This coherence requirement saves us from guessing which hypothesis is true when designing optimal experiments.

We can go further and consider what specific evidence functions satisfying (\ref{eqn:observe_coherence}) are particularly appealing. We want our evidence function to be probability based because hypothesis testing is fundamentally about seeking probabilistic evidence, usually in the form of {\it p-}values or posterior probabilities. Indeed, for the purposes of test information, the traditional estimation information focus on variance and spread is in general inadequate. From the Bayesian perspective, a sensible probability based evidence function is
\begin{align}
\declr(z;H_0,H_1)=\frac{z}{\pi_1+\pi_0z},
\label{eqn:postef_h0}
 \end{align}
where $\pi_0$ and $\pi_1$ are the prior probabilities of $H_0$ and $H_1$, respectively (for simplicity we assume $\pi_0+\pi_1=1$). When $z = \bfac(\datx|H_0,H_1)$,  (\ref{eqn:postef_h0}) is just the posterior to prior probability ratio for $H_0$. The resulting dual expected test information measures are
\begin{align}
\bdecinfo(\ept;H_0,H_1) &= 1 - \Ex\left[\left.\frac{\zlr(\ranx)}{\pi_1+\pi_0\zlr(\ranx)}\right|H_1\right] \label{eqn:post1}\\
\bdecinfo(\ept;H_1,H_0) &= 1 - \Ex\left[\left.\frac{1}{\pi_1+\pi_0\zlr(\ranx)}\right|H_0\right], \label{eqn:post0}
\end{align}
where $\zlr(\ranx)$ denotes the Bayes factor $\bfac(\ranx|H_0,H_1)$. The measure (\ref{eqn:post1}) is simply the expected difference between the prior and posterior probability of the null, relative to the prior probability, when the data are actually from the alternative. That is, (\ref{eqn:post1}) constitutes the relative loss in probability of the null. The measure (\ref{eqn:post0}) is the same but with the roles of $H_0$ and $H_1$ switched. Since (\ref{eqn:observe_coherence}) is satisfied, the coherence identity (\ref{eqn:info_link}) tells us that $\bdecinfo(\ept;H_0,H_1)=\bdecinfo(\ept;H_1,H_0)$. This and the straightforward Bayesian probability interpretation of (\ref{eqn:post1}) make (\ref{eqn:postef_h0}) a particularly appealing choice of evidence function.

There are also other examples of evidence functions that satisfy (\ref{eqn:observe_coherence}), e.g.,
\begin{align}
\declr(\zlr) = \frac{1}{2}\log(\zlr) - \frac{1}{2}\zlr\log(\zlr).
\label{eqn:lambdadiv0}
\end{align}
For this evidence function, both sides of (\ref{eqn:info_link}) equal $\frac{1}{2}KL(f(\cdot|H_1)||f(\cdot|H_0))+ \frac{1}{2}KL(f(\cdot|H_0)||f(\cdot|H_1))$. Historically, this symmetrized form of KL divergence is the divergence that \cite{kl1951} originally suggested (without scaling by a half). 
Intuitively, it can be interpreted as a measure of the expected test information when the two hypotheses are considered equally likely apriori. However, symmetrized KL divergence does not have a straightforward probability interpretation, and therefore we prefer (\ref{eqn:postef_h0})-(\ref{eqn:post0}). 

\subsection{Nuisance parameters}
\label{sec:nuisance}

Many statistical problems come with nuisance parameters. In the frequentist setting, once data have been observed, estimates of the nuisance parameters can be inserted to give a point estimate of the expected test information (\ref{eqn:freq_decinfo}). A confidence interval for (\ref{eqn:freq_decinfo}) can be obtained by evaluating it for values of the nuisance parameters within a confidence interval. (Both could instead be done for observed test information, see Section \ref{sec:observed_measures}.) In design problems, data are typically yet to be collected but (\ref{eqn:freq_decinfo}) could be evaluated on a grid of values of the nuisance parameters.

In the Bayesian context, the nuisance parameters $\nuis_0$ (under the null) and $\nuis_1$ (under the alternative) are simply integrated out along with the parameters that define the hypotheses. That is,
\begin{align}
\bdinfofull(\ept;H_0,H_1,\prior,\nuispriornull,\nuisprioralt) = \decz - \Ex & \left[\left.\declr\left(\bfac(\ranx|H_0,H_1)\right)\right|H_1\right], \label{eqn:berger_approach}
\end{align}
where the Bayes factor is now given by
\begin{align}
\frac{\int_\psupnull\int_{\nsupnull} f(\ranx|\para,\nuisnull)\nuispriornull(\nuisnull|\para)\prior(\para|H_0)d\nuisnull d\para}{\int_\psupalt\int_{\nsupalt} f(\ranx|\para,\nuisalt)\nuisprioralt(\nuisalt|\para)\prior(\para|H_1)d\nuisalt d\para},\label{eqn:nuisance_bf}
\end{align}
with $\nsup_i$ being the support of the prior density $\nuisprior_i$ of $\nuis_i$, for $i=0,1$. Clearly, the mathematical properties of (\ref{eqn:berger_approach}) are the same as those of (\ref{eqn:bayes_decinfo}). 

As mentioned in Section \ref{sec:expected_test}, alternatives to the Bayes factor in Definition \ref{def:decinfo} can be used at the expense of the coherence of the Bayesian method and simplicity. For example, those intending to use the likelihood ratio test, may opt to mimic the likelihood ratio test statistic by calculating
\begin{align}
\bdinfofull(\ept;H_0,H_1,\prior,\nuisprioralt) = \decz - \Ex\left[\left.\declr\left(\frac{f(\ranx|\para\mlenull,\nuismlenull)}{f(\ranx|\para\mlealt,\nuismlealt)}\right)\right| H_1\right], \label{eqn:max_like_approach}
\end{align}
where $\para\mlei$ and $\nuismlei$ are the MLEs of $\para$ and $\nuis_i$, respectively, under hypothesis $H_i$, for $i=0,1$. 
In this work we focus on the expected test information given in Definition \ref{def:decinfo} (and (\ref{eqn:freq_decinfo})) and thus leave the theoretical investigation of measures such as (\ref{eqn:max_like_approach}) for future work. However, we include numerical results  based on (\ref{eqn:max_like_approach}) in Section \ref{sec:binary}.

\subsection{Probit and complementary log-log regression example}
\label{sec:binary}

Consider the binary regression model
\begin{align}
\ranx_1, \ldots, \ranx_n\ |M,\nuis_\para, g_\para \overset{i.i.d}{\sim} \bern(p_i),\label{eqn:logistic_glm}
\end{align}
where
\begin{align}
M^T=\left(\begin{matrix} 1 & 1 & \cdots & 1 \\
t_1 & t_2 & \cdots & t_n\end{matrix}\right)
\label{eqn:simple_design}
\end{align}
is the design matrix (i.e., essentially $\ept$),  and $g_\para(p_i) = \nuisinttheta + \nuisslopetheta t_i$, for the link function $g_\para$ and $\para \in \{0,1\}$. The sharp hypotheses of interest are $H_0:\para=0$ and $H_1:\para=1$, where $g_0(p) = \log(-\log(1-p))$ and $g_1(p)=\Phi^{-1}(p)$  are the familiar complementary log-log and Probit link functions, respectively ($\Phi$ is the standard Normal cumulative distribution function). In this model selection scenario, the coefficients $\nuisnull=\{\nuisnullint,\nuisnullslope\}$ and $\nuisalt=\{\nuisaltint,\nuisaltslope\}$ are nuisance parameters, and we assign the prior distribution
\begin{align}
\nuis_i|H_i &\sim N_2(\altmean,\altcov),\label{eqn:nuis_prior}
\end{align}
for $i=0,1$, with $\eta$ and $R$ to be determined. Our independent prior for $\para$ is Bernoulli$(1/2)$. The design problem is to choose the design which will provide the most information for distinguishing between the two link functions. \citet{ponce_thesis} previously considered a similar problem but formulated it as an estimation design problem by parameterizing a continuum of models between the logistic and complementary log-log binary regression models, as for example in \citet{glm}, page 378. 

We consider designs of 5 unique points in $[-1,1]$ with 100 replications of each. Within this class, we optimize the expected test information under the \postv and \logv evidence functions, i.e., (\ref{eqn:postef_h0}) and $\declr(\zlr) = \log(\zlr)$. Since $\nuisnull$ and $\nuisalt$ are nuisance parameters, we use the two measures in Section \ref{sec:nuisance}, i.e., the Bayes and MLE plug-in approaches. Under the \postv evidence function, we denote these two measures  by $\critpb(M)$ and $\critpl(M)$, respectively. Similarly, under the \logv evidence function the two measures are denoted $\crittb(M)$ and $\crittl(M)$, respectively. The $\crit$ notation indicates a design criterion, $P$ indicates a connection to the expected posterior probability of $H_1$, given by $\pi_1 + \pi_0 \critpb(M)$, and $TK$ indicates the testing context and the KL divergence between the marginal data models. The criteria are computed using Monte Carlo simulation and the optimal design under each criterion is found using a single point exchange algorithm similar to that introduced by \cite{fedorov1972}. The design matrix optimizing criterion $C$ is denoted $M_{C}$, for $C\in\{\postbayes,\postlrt,\tkbayes,\tklrt\}$, collectively called the \poptw \space and \tkw \space designs.

\begin{figure}[t]
\includegraphics[width=\textwidth,trim = 0mm 12.5mm 0mm 5mm]{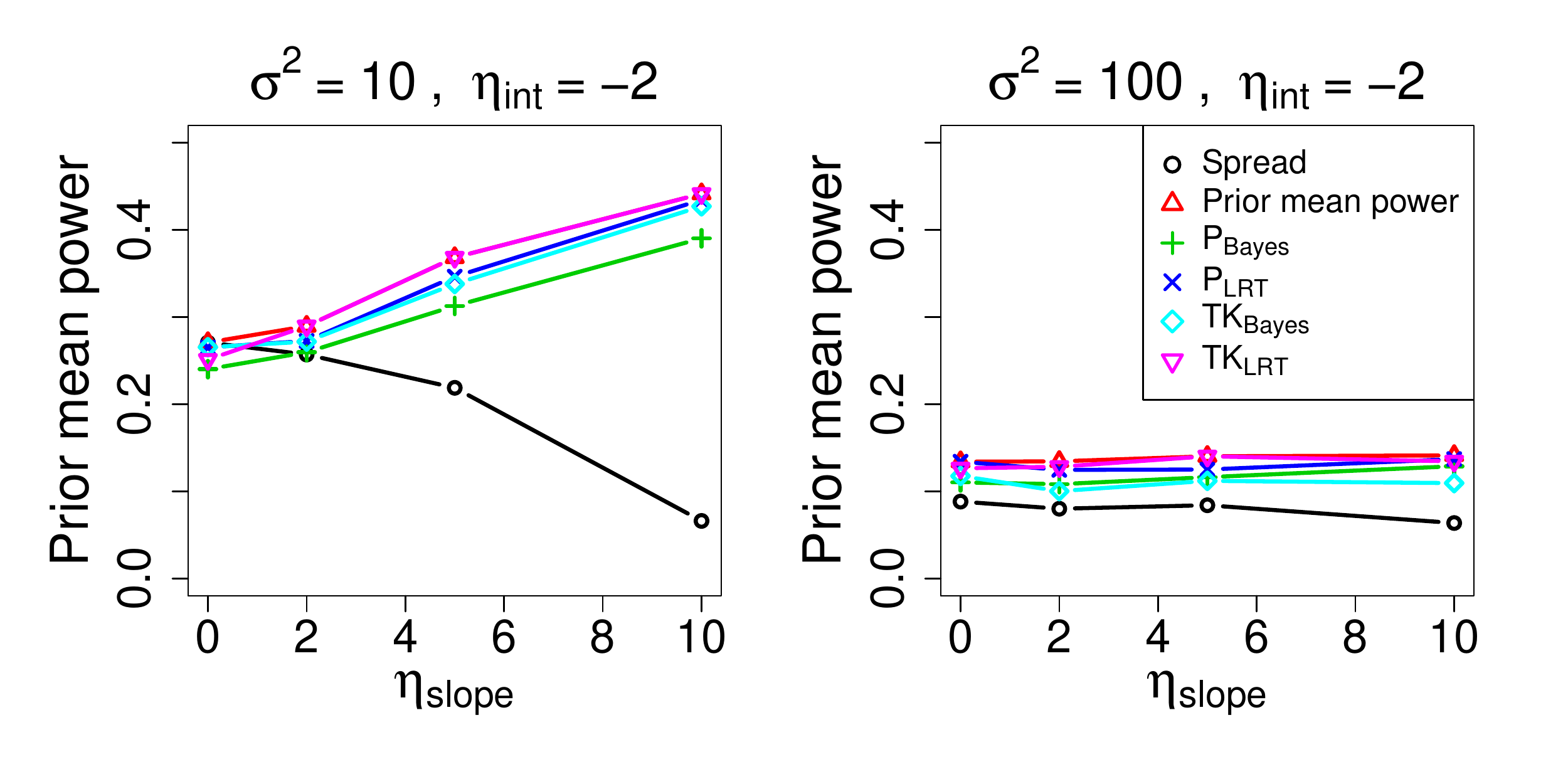}
\caption{Prior mean power of the likelihood ratio test under $M_{C}$, for $C\in\{\text{Spread},\text{Power},$ $\postbayes,\postlrt,\tkbayes,\tklrt\}$, for different settings of the priors in (\ref{eqn:nuis_prior}). For the \poptw \space designs we set $\pi_0=\pi_1=0.5$.\label{fig:binary_power}}
\end{figure}

We need a separate criterion by which we can evaluate and compare the optimal designs. Since power is a common quantity of interest, we choose the criterion to be the prior mean power of the likelihood ratio test, i.e.,
\begin{align}
\int_{\Omega} \mathcal{P}(M;\nuisnull,\nuisalt)\nuispriornull(\nuisnull|H_0)\nuisprioralt(\nuisalt|H_1)d(\nuisnull,\nuisalt),
\label{eqn:prior_mean_power}
\end{align}
where $\Omega=\nsupnull \times \nsupalt$, and $\mathcal{P}(M;\nuisnull,\nuisalt)$ denotes the power of the likelihood ratio test under design matrix $M$ and given the parameters $\nuisnull$ and $\nuisalt$ (for a test size of 5\%). Section \ref{sec:power} discusses reasons why (\ref{eqn:prior_mean_power}) or similar summaries of power are not the only measures of expected test information, or even particularly good measures. Nonetheless, the relative familiarity of (\ref{eqn:prior_mean_power}) makes it suitable for our current purpose of comparing the performance of the different optimal designs. Figure \ref{fig:binary_power} shows the prior mean power under $M_{C}$, for $C\in\{\postbayes,\postlrt,\tkbayes,\tklrt\}$, given various specifications of the priors in (\ref{eqn:nuis_prior}). In all cases the covariance matrix $\altcov=\sigma^2I_2$ and only $\sigma^2$ and $\altmean$ are indicated. Note that we tried several values of $\altmeanint$ but the results were qualitatively very similar, so Figure \ref{fig:binary_power} only shows results for $\altmeanint=-2$.

Also shown is the prior mean power under $\mpower$ and $\mspread$, the maximum prior mean power design and the spread of points $-1,-0.5,0,0.5,1$ (replicated 100 times), respectively. The \poptw \space and \tkw \space designs all perform well in terms of prior mean power, and in some cases yield considerably greater prior mean power than $\mspread$.  For example, when $\sigma^2=10$ and $\altmean=\{-2,10\}$ (left panel of Figure \ref{fig:binary_power}), the design $\mspread$ has prior mean power $0.07$ while the \poptw \space and \tkw \space designs are all relatively close to achieving the maximum prior mean power of $0.44$. The problem with $\mspread$ in this case is that both inverse link functions practically go from 0 to 1  over a small range of the covariate $t$ and therefore spreading the design points over the whole interval $[-1,1]$ is not an effective strategy.

\begin{figure}[t]
\includegraphics[width=0.95\textwidth,trim = 10mm 2.5mm 5mm 0mm]{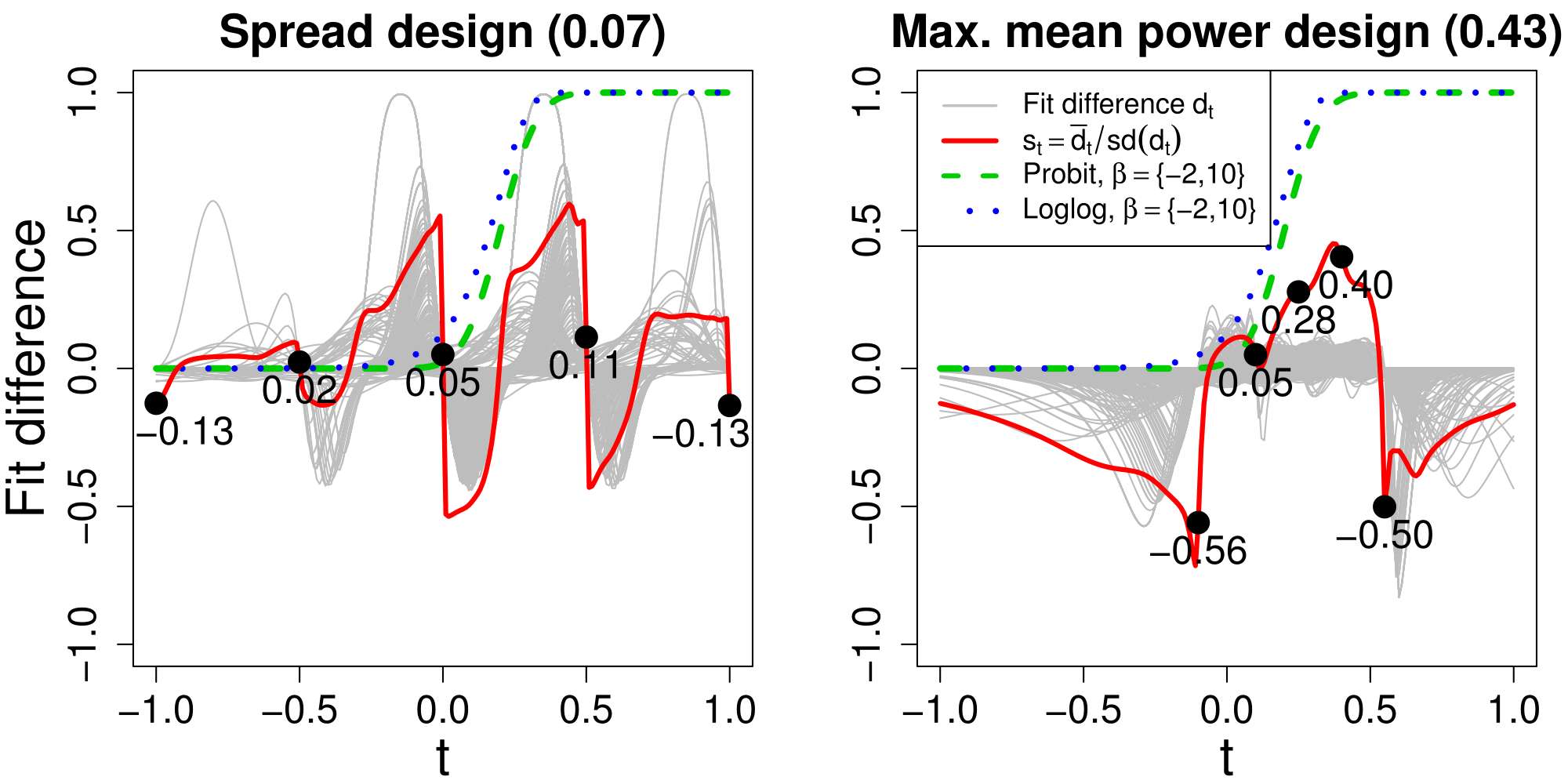}
\caption{Comparison of $\mspread$ (left) and $\mpower$ (right), with the parameters of (\ref{eqn:nuis_prior}) set to $\altmean=\{-2,10\}$ and $R=10I_2$. Thin grey lines show $d_t$, for each simulated dataset, and the thick red line shows $s_t$ (both are described in the main text). The large dots show the design point locations (x-coordinates) and the corresponding values of $s_t$ (y-coordinates and numbers below).  \label{fig:binary_design_comparison}}
\end{figure}

To investigate the designs further, Figure \ref{fig:binary_design_comparison} compares $\mspread$ and $\mpower$. In each plot,  the design is given by the x-coordinates of the large dots (100 binary observations are recorded at each). In this illustration, 500 datasets were simulated under each design with $\altmean=\{-2,10\}$ and $\sigma^2=10$. For reference, the Probit (dotted blue line) and complementary log-log (dashed green line) inverse link functions are plotted for $\nuisnull=\nuisalt=\{-2,10\}$. However, since $\sigma^2>0$, the actual value of $\nuisnull$ and $\nuisalt$ vary across the simulated datasets. Furthermore, for any given dataset, there is uncertainty associated with the MLE of $\nuisnull$ and $\nuisalt$. These two sources of variation are captured by the spread of the solid thin grey lines in Figure \ref{fig:binary_design_comparison}; each corresponds to a single simulated dataset and traces the fit difference $d_t(x^{(j)}) = g_1^{-1}(\nuisaltinthatj+\nuisaltslopehatj t)-g_0^{-1}(\nuisnullinthatj+\nuisnullslopehatj t)$ for $t \in [-1,1]$, where $\nuisinthatji,\nuisslopehatji$ are the MLEs of $\nuisinti,\nuisslopei$ for dataset $x^{(j)}$, for $i=0,1,$ and $j=1,\dots,500$. The distribution of the fit differences at any point $t$ indicates our ability to distinguish the two inverse link functions at that point based on maximum likelihood fits. The solid thick red line summarizes by tracing the relative mean fit difference $s_t = \bar{d}_t/\mbox{sd}(d_t)$, where $\bar{d}_t$ and $\mbox{sd}(d_t)$ are the mean and standard deviation of $d_t$ over 500 simulations, respectively.

The y-coordinates of the large dots give the values of $s_t$ at the design points (as do the numbers below the large dots). 
As expected, $|d_t|$ is generally small at the design points, but for $\mpower$ the variability in $d_t$ is low and thus $|s_t|$ is larger at the design points than under $\mspread$. The low variability is achieved by grouping the design points together near the important steep section of the reference inverse link functions. The complementary log-log  and Probit regression models fit by maximum likelihood are known to differ principally in the tails, and hence $s_t$ is not largest at the design points in the central steep section. However, these points constrain the fits, thus reducing variability in $d_t$ so that the two design points in the tails have large values of $s_t$. 
The designs $M_{\tkbayes}$ and $M_{\tklrt}$ are almost identical to $\mpower$, which is to be expected because intuitively the prior mean power should increase as the expected negative log Bayes factor (or likelihood ratio) increases.

\subsection{Normal linear regression coefficient tests}
\label{sec:testb}

We now discuss the Normal linear regression model
\begin{eqnarray}
\ranx|\nuis,M &\sim& N_n(M\nuis,\sigma^2I),\label{eqn:linearb}
\end{eqnarray}
and the hypotheses $H_0:\nuis=\nuisnull$ and $H_1:\nuis \sim N_d(\altmean,\sigma^2\altcov)$. The goal is to test the adequacy of a specific value $\nuisnull$ of the regression coefficients, rather than treating them as nuisance parameters as we did in Section \ref{sec:binary}.  Thus, $\nuis$ is now playing the role of $\para$ in the expected test information of Definition \ref{def:decinfo}, and formally we restrict its support under $H_1$ to be $\nsalt = \Bbb{R}^d/\{\nuisnull\}$. We again consider the criteria $\critpb(M)$ and $\crittb(M)$ which, since there are no nuisance parameters, are now simply given by Definition \ref{def:decinfo} under the \postv and \logv evidence function, respectively.

In the linear regression setting, $\crittb(M)$ has the closed form
\begin{align}
\label{eqn:bayes1_utility}\crittb(M;\nuisnull,\altmean,\altcov) =  \int_\xsup\log\left(\frac{f(x|H_1,M)}{f(x|\nuisnull,M)}\right)f(x|H_1,M)dx \hspace{1cm}\\
 \nonumber= \frac{1}{2}\left(\frac{1}{\sigma^2}(\altmean-\nuisnull)^TM^TM(\altmean-\nuisnull) + \mbox{tr}( M^TM\altcov)- \log\left(|I+M\altcov M^T|\right)\right).
\end{align}
The first term of (\ref{eqn:bayes1_utility}) confirms our intuition that the expected test information is large when $\nuisnull$ and the mean of the alternative are well separated (with respect to $(M^TM)^{-1}$). Heuristically, the second term of (\ref{eqn:bayes1_utility}) tells us to maximize the ``ratios" of the prior  (alternative) variance of each parameter to the regression estimate variance. This is intuitive because we need the estimation variance to be small in comparison to the prior variance in order to effectively distinguish $\nuis$ and $\altmean$ (and hence further distinguish $\nuis$ from $\nuisnull$).  The final term penalizes the alternative for introducing uncertainty in $\nuis$, i.e., for avoiding exclusion of the true model by being vague.

The \tkc \space (\ref{eqn:bayes1_utility}) is closely related to the D-optimality criterion popular in estimation problems (e.g., see the review by \citealt{chaloner}). The D-optimality criterion is derived from the expected estimation information suggested by \cite{lindley1956}, and is given by
\begin{eqnarray}
\crit_D(M)=-\log|V|, \label{eqn:dopt}
 \end{eqnarray}
where $V = \sigma^2(M^TM+R^{-1})^{-1}$ is the posterior covariance matrix of $\nuis$ (for any value of $\ranx$). The criteria $\crit_D$ and $\crittb$ are both entropy based, but the dependence of $\crittb$ on $\nuisnull$ and $\altmean$ distinguishes this criterion from $\crit_D$ and other estimation focused criteria.

To gain some intuition, let us consider a simple linear regression with
 \begin{eqnarray}
M^T=\left(\begin{matrix} 1 & 1 & \cdots & 1 \\
t_1 & t_2 & \cdots & t_{n}\end{matrix}\right), \quad R=\left(\begin{matrix} \sigma^2_{\text{int}} & \covoff  \\
 \covoff & \sigma^2_{\text{slope}}\end{matrix}\right).
\end{eqnarray}
In this scenario, it is well known that the \doptw \space design places half of the design points $t_i$, $i=1,\dots,n$, at $1$ and the other half at $-1$ (or, if $n$ is odd, $(n+1)/2$ points at one boundary and $(n-1)/2$ at the other). Let $\Delta = (\altmeanint-\nuisnullint)(\altmeanslope-\nuisnullslope) + \covoff$, where $\altmeanint$ and $\altmeanslope$ are the mean intercept and mean slope of the alternative model, respectively. If $\covoff=0$, then the sign of $\Delta$ tells us if the lines $\nuisnullint+\nuisnullslope t$ and $\altmeanint+\altmeanslope t$ have greater separation at $-1$ or at $1$. 
For any $\covoff$, it is easily shown that $\crittb$ is optimized by placing  all points at $1$ if $\Delta>0$, by placing them at $-1$ if $\Delta<0$, and by any design dividing the points between the boundaries if $\Delta=0$. 
Generally, designs based on test information measures trade robustness for power in distinguishing particular models, and the behavior just described is an instance of the inevitable sensitivity to the hypotheses mentioned in Section \ref{sec:intro}. However, in the current context, we found that designs optimizing $\critpb$ are slightly more robust than those optimizing $\crittb$ in that they divide the points between both the boundaries, unless the hypotheses are far more separated at one boundary than at the other. 

\section{Observed Test Information in Theory and Applications}
\label{sec:observed}

\subsection{Observed test information: building blocks}
\label{sec:observed_measures}

Observed test information is key in practice when we observe some data and wish to know how much information they contain in order to decide if we should collect more. It is also important conceptually because it is the implicit building block for expected  and conditional test information.

Consider the setting in Section \ref{sec:degroot}. After an experiment $\ept$ is conducted, the observed estimation information gained is  the reduction in uncertainty, $\uncp(\prior) - \uncp(p(\cdot|\datx))$, where $\datx \in \xsup$ is the observed outcome.  Observed estimation information is not necessarily  non-negative because, by chance, after observing $\datx$ we may have more uncertainty  about $\para$ as measured by $\uncp$, e.g., the posterior may be more diffuse than the prior due to likelihood-prior conflict; see \cite{object_prior}. (This posterior ``dilation" can even be deterministic; see \cite{wasserman1993}.) 
\cite{degroot1962} did not explicitly mention observed information, but \cite{lindley1956} did define it (as above) when $\uncp$ is the entropy function (\ref{eqn:entropy}). \citet{ginebra2007} restricted all observed information measures to be non-negative and asserted that this allows them to be interpreted as capturing model checking information, in addition to information about $\para$. However, \citet{ginebra2007} did not explain why adding a non-negativity condition is necessary or sufficient to ensure that a model checking interpretation is reasonable.  From a Bayesian perspective, the definition given by \cite{lindley1956} is valid, and we therefore take this as the basis of observed estimation information.

Following \cite{lindley1956}, we define observed test information in Definition \ref{def:obs_decinfo} (below) by simply removing the expectation appearing in the expected test  information of Definition \ref{def:decinfo}. However, the resulting relationship between observed and expected information is more subtle than in the estimation case. Indeed, Definition \ref{def:decinfo} conditions on $H_1$ to average over the unobserved data, but the actual data used in Definition \ref{def:obs_decinfo} may be generated under $H_0$. 

\begin{defi}
\label{def:obs_decinfo}
The observed test information provided by the experiment $\ept$ for comparing the hypotheses $H_0: \para \in \psupnull$ and $H_1: \para \in \psupalt$, for a given evidence function $\declr$ and a proper prior $\prior$, is defined as
\begin{eqnarray}
\label{eqn:obs_bdecinfo}
\bdinfofull(\ept;H_0,H_1,\prior,\datx) = \decz - \declr(\bfac(\datx|H_0,H_1)),
\end{eqnarray}
where $\datx$ is the observed outcome of $\ept$, and $\psupnull \cap \psupalt = \emptyset$.
\end{defi}

Since Bayesians condition on observed data, the prefix `observed' is redundant, but it is retained for clarity. The quantity defined by (\ref{eqn:obs_bdecinfo}) is not necessarily non-negative. However, it is positive when $\declr$ is increasing and the Bayes factor favors $H_1$, i.e., $\bfac(\datx|H_0,H_1)<1$. Often, it seems sensible for $\declr$ to be increasing because we want (\ref{eqn:obs_bdecinfo}) to increase as the Bayes factor decreases towards zero (since observed test information should be compatible with Definition  \ref{def:decinfo} which assumes $H_1$).  For $\declr$ increasing, a negative value of $\bdinfofull(\ept;H_0,H_1,\prior,\datx)$ indicates that the evidence in the observed data supports $H_0$, either because $H_0$ is in fact the more accurate hypothesis or due to chance. Since the data can support only one of the hypotheses, for increasing $\declr$ it follows that exactly one of the dual observed test information measures $\bdinfofull(\ept;H_0,H_1,\prior,\datx)$ and $\bdinfofull(\ept;H_1,H_0,\prior,\datx)$ will be positive (unless they are both zero). Also, usually only one of $\bdinfofull(\ept;H_0,H_1,\prior,\datx)$ and $\bdinfofull(\ept;H_1,H_0,\prior,\datx)$ will reasonably approximate the corresponding expected test information, $\bdinfofull(\ept;H_0,H_1,\prior)$ and $\bdinfofull(\ept;H_1,H_0,\prior)$, respectively. 


We highlight that, in the current observed data case, our use of dual measures is again key because it ensures a symmetric treatment of the hypotheses, which is not easily achieved by other means. For example, consider a $\declr$ that is concave, increasing, and passes through $\{0,1\}$ (for all non-zero prior probabilities $\pi_0$ and $\pi_1$ whose sum is one), then $\declr(\zlr) + \zlr \declr(\zlr^{-1})$ is also concave, it gives expected test information  $\bdinfofull(\ept;H_0,H_1,\prior)+\bdinfofull(\ept;H_1,H_0,\prior)$, and yields non-negative observed test information, as required by \citet{ginebra2007}. However, in many cases, excluding the case of (\ref{eqn:lambdadiv0}), this approach does not treat the hypotheses equally. For example, we can modify the evidence function $\declr(\zlr)=\sqrt{\zlr}-1$ to $\declr(\zlr) + \zlr \declr(\zlr^{-1})=2\sqrt{\zlr}-1-\zlr$, but the resulting observed test information has a maximum of one for data supporting $H_1$, and is unbounded for data supporting $H_0$. Our approach using dual observed test information measures is therefore more appealing because $\bdinfofull(\ept;H_0,H_1,\prior,\datx)$ and $\bdinfofull(\ept;H_1,H_0,\prior,\datx)$ are symmetrically defined. 

Next, in the same spirit, we define the observed conditional test information from experiment $\ept_2$ after observing  $\datx_1$ of experiment $\ept_1$ to be
\begin{eqnarray}
\label{eqn:ob_con_decinfo}
\bdinfofull(\ept_2|\ept_1;H_0,H_1,\prior,\datx_1)= \dec(\datx_1) - \Extwo[\dec(\datx_1,\ranx_2)|H_1,\datx_1].
\end{eqnarray}
This is simply the information given in Definition \ref{def:con_decinfo}, but without an expectation over $\datx_1$. In sequential design we require a version of the coherence identity (\ref{eqn:info_link}) to hold for (\ref{eqn:ob_con_decinfo}). In particular, given some observed data $\datx_1$, we want the optimality of our design to be free of
the validity of hypothesis $H_0$ or $H_1$. The symmetry condition (\ref{eqn:observe_coherence}) implies that
\begin{eqnarray}
\decinfo(\ept_2|\ept_1;H_0,H_1,\prior,\datx_1)= \zlr(\datx_1)\decinfo(\ept_2|\ept_1;H_1,H_0,\prior,\datx_1),
\label{eqn:con_coherence}
\end{eqnarray}
for all $\datx_1 \in \xsup$, where $\zlr(\datx_1)=\bfac(\datx_1|H_0,H_1)$. The factor $\zlr(\datx_1)$ appears in (\ref{eqn:con_coherence}) (but not in (\ref{eqn:observe_coherence})) because the observed data $\datx_1$ already favors one of the hypotheses before any new data are collected. If (\ref{eqn:con_coherence}) holds, then a design that optimizes (\ref{eqn:obs_bdecinfo}) also optimizes $\decinfo(\ept_2|\ept_1;H_1,H_0,\prior,\datx_1)$. Hence we do not need to know which hypothesis is true when choosing $\ept_2$.

\subsection{Sequential design for linear regression coefficient tests}
\label{sec:seq}

\begin{figure}[t]
\includegraphics[width=0.9\textwidth,trim=10mm 5mm 10mm 5mm]{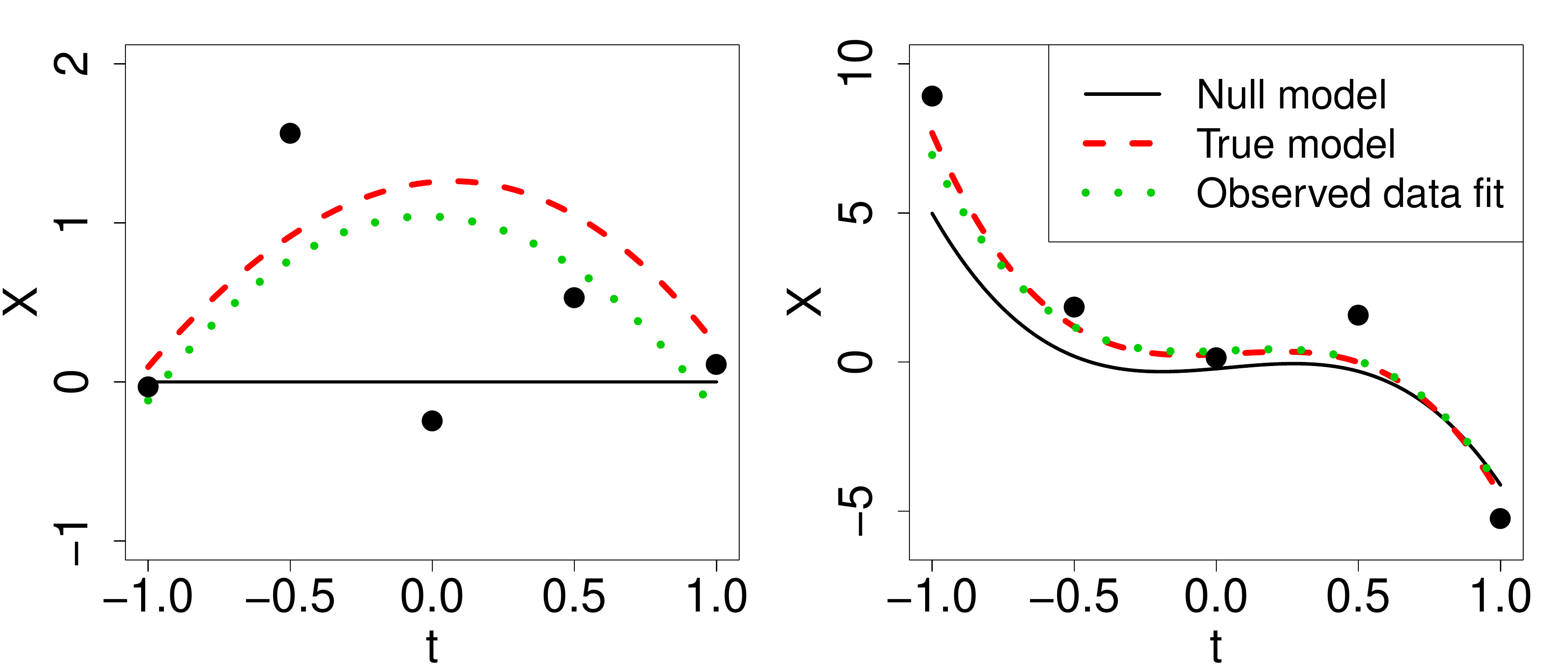}
\caption{The null and true  cubic regression models and the observed data posterior mean fit. The observed data are indicated by large dots. The left and right plots show example simulations used in producing parts (a) and (b) of Figure \ref{fig:opt_point_spread_results}, respectively. \label{fig:cubic_model}}
  \end{figure}

Consider the linear regression model (\ref{eqn:linearb}) introduced in Section \ref{sec:testb} and the test of $H_0: \nuis=\nuisnull$ against $H_1: \nuis\sim N(\eta,R)$ (i.e., $\sigma^2=1$). Given some initial observed data $\datx\ob$, the sequential design problem is to choose a design matrix $M\mis$ for additional data $\ranx\mis$. In our simulation study, we generate a parameter vector $\nuis$ under $H_1$ and then simulate the initial observed data $\datx\ob$ according to a cubic regression model of the form (\ref{eqn:linearb}), i.e., the design matrix is
 \begin{eqnarray}
M\ob^T=\left(\begin{matrix} 1 & 1 & \cdots & 1 \\
t_1 & t_2 & \cdots & t_{n\ob} \\
t_1^2 & t_2^2 & \cdots & t_{n\ob}^2\\
t_1^3 & t_2^3 & \cdots & t_{n\ob}^3\end{matrix}\right),
\label{eqn:cubic_design}
\end{eqnarray}
 where $t_i\in[-1,1]$ are the design points, for $i=1,\dots,n\ob$. Specifically, we set $n\ob=5$ and the observed data design points $t_i$, $i=1,\dots,5$, are $-1$,$-0.5$, $0$, $0.5$, $1$. Examples of $\datx\ob$ are plotted in Figure \ref{fig:cubic_model}.

Given the observed data, $n\mis=5$ new design points are chosen by optimizing the observed conditional test information (\ref{eqn:ob_con_decinfo}) with respect to the design matrix $M\mis$, under the \postv and log evidence functions. That is, we optimize the conditional versions of the \popt \space and \tkcs \space discussed in Section \ref{sec:testb}. The conditional \poptc \space is the expected reduction in posterior probability of the null when we collect $\ranx\mis$ relative to its prior probability. We approximate it using a Monte Carlo estimate, under the prior probabilities $\pi_0=\pi_1=0.5$. The conditional \tkc \space is straightforwardly given by $\log(\zlr(\datx\ob)) + \crittb(M;\nuisnull,\altmean\ob,V\ob)$,
where $\crittb$ is specified in (\ref{eqn:bayes1_utility}), and $\altmean\ob$ and $V\ob$
are the observed data posterior mean and covariance matrix of $\nuis$ under $H_1$, respectively. For comparison, we also optimize the (deterministic) conditional D-optimality criterion $\log|V\ob|+\log\left|\left(M\mis\right)^TM\mis+V\ob^{-1}\right|$ with respect to $M\mis$. Often it is not clear how to use the D-optimality criterion and other estimation based criteria to choose designs for testing, but the current scenario is an exception because the hypotheses are nested.

\begin{figure}[t]
\includegraphics[width=0.9\textwidth,trim= 0mm 10mm 0mm 0mm]{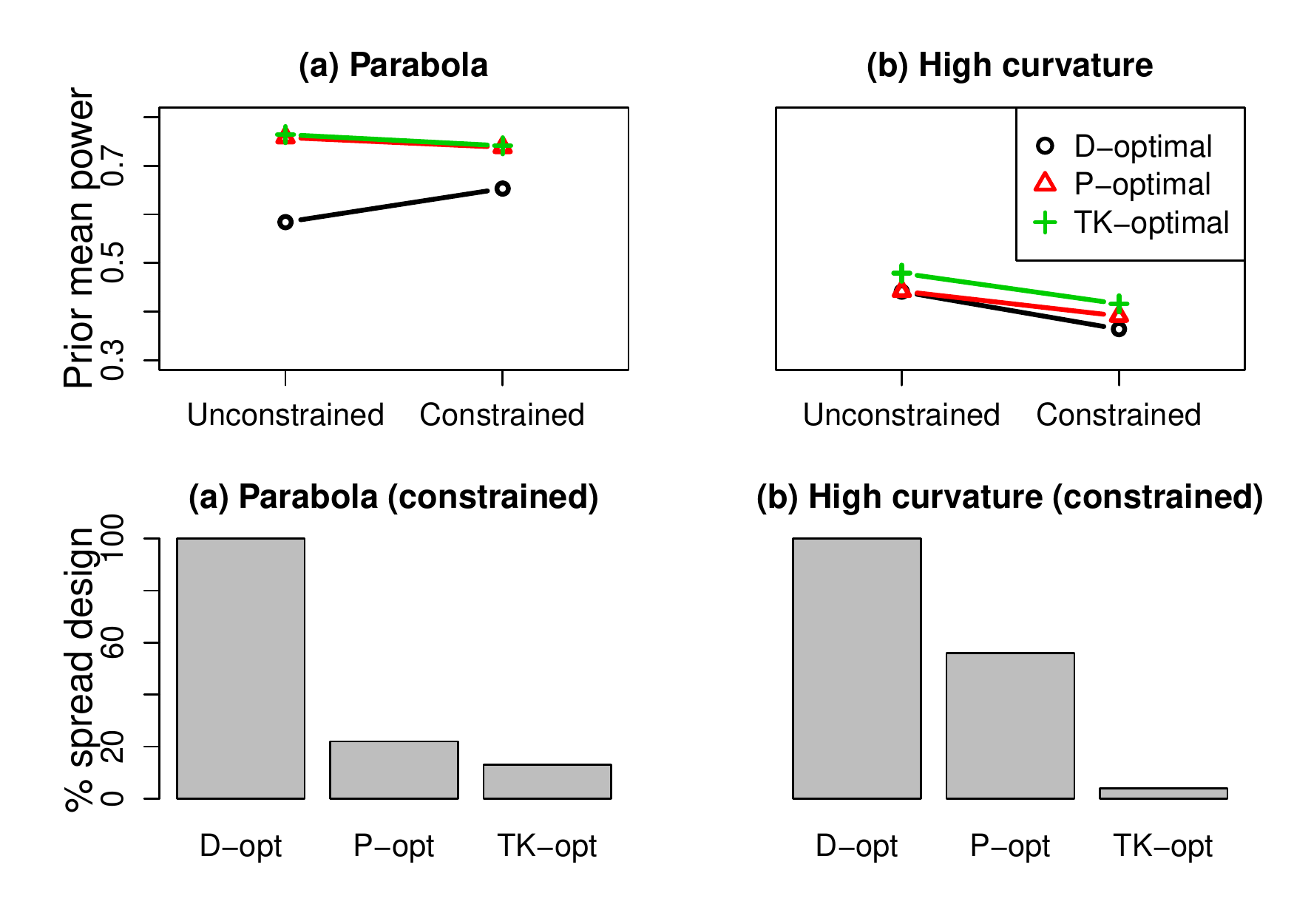}
\caption{
Prior mean power of the likelihood ratio test under the conditional \dopt, \popt, and \tk \space procedures, across 250 datasets  simulated under $\nuis \sim N(\eta,0.2I_4)$ (first row, unconstrained values). The main text describes the generation of $\nuisnull$ and $\altmean$ for parts (a) and (b). 
The first row constrained values show the prior mean powers when the only missing data designs allowed are (i) and (ii) (see the main text). For this case, the second row shows the percentage of simulations in which design (i) was selected. \label{fig:opt_point_spread_results}}
\end{figure}

To generalize beyond a single value of $\nuis$, we generate $\nuis^{(j)}\sim N(\eta,R)$, for $j=1,\dots,100$, and for each $j$ we generate observed datasets $\datx\ob^{(j,k)}$, for $k=1,\dots,250$. Then, for each simulated dataset $\datx\ob^{(j,k)}$, we find the conditional \poptw, \tkw, and \doptw \space design for the missing data $\ranx\mis$. To compare performance, we also calculate the prior mean power (\ref{eqn:prior_mean_power}) of the likelihood ratio test under each of these three procedures. 
In our simulations, we set $R=0.2I_4$, and use various values of $\eta$ and $\nuisnull=(\nuisnullint,\nuisnulllin,\nuisnullquad,\nuisnullcubic)$. First, Figure \ref{fig:opt_point_spread_results} part (a) corresponds to simulations with $\nuisnull = (0,0,0,0)$ and $\eta = (1.1,0,-1.3,0)$, i.e., the alternative mean model is parabola shaped. For these choices, the maximum separation between the null and true model is usually not at the boundaries of the interval $[-1,1]$; see the example simulation given on the left of Figure \ref{fig:cubic_model}. The top row of Figure \ref{fig:opt_point_spread_results} shows the prior mean power of the three procedures when all design points in $[-1,1]$ are allowed (i.e., $M\mis$ is unconstrained) and also when only two possibilities for $M\mis$ are allowed (these latter results are for the constrained optimization example discussed shortly). For part (a), the conditional \doptw \space procedure performs relatively poorly because it places all the new points at the boundaries, a good strategy for estimation but not for the current hypothesis test. The conditional \popt \space and \tk \space procedures instead place the points near $t=0$, and consequently are substantially superior in terms of prior mean power.

For the simulations corresponding to Figure \ref{fig:opt_point_spread_results} part (b), we generated $\nuisnullint^{(j)}, \nuisnulllin^{(j)} \sim \text{Uniform}(-1,1)$ and $\nuisnullquad^{(j)},\nuisnullcubic^{(j)}\sim\text{Uniform}(-10,10)$, and then set $\altmean^{(j)}=\nuisnull^{(j)}$ and drew $\nuis^{(j)}\sim N(\altmean^{(j)},R)$, for $j=1\dots 100$. Under these settings, the maximum separation between the null curve and the observed data posterior mean fit tends to be at one of the boundaries, and therefore the conditional \dopt \space procedure performs reasonably. 
Thus, in part (b) of Figure \ref{fig:opt_point_spread_results} the three procedures perform similarly. 

We now briefly investigate how the three criteria perform if we impose some robustness to model misspecification. The points labeled ``constrained" in the first row of Figure \ref{fig:opt_point_spread_results} show the prior mean power of the likelihood ratio test when the three criteria are used to choose between two missing data designs: (i) the spread of points $\boldsymbol{t}_{\text{\tiny spread}}=\{-1,-0.5,0,0.5,1\}$; (ii) the narrower spread of points $\frac{1}{5}\boldsymbol{t}_{\text{\tiny spread}}+\text{sep}_{\max}$, where $\text{sep}_{\max}$ is the location of maximum separation between the observed data posterior mean model and the null mode. If $\text{sep}_{\max}$ is near a boundary then all the points are shifted left or right to avoid any crossing the boundary, but they still cover an interval of length 0.4. The results follow a similar pattern to before, except that now the prior mean power is usually lower, principally because designs placing all the points at a single location have been excluded. The first row of Figure \ref{fig:opt_point_spread_results} shows that the constrained conditional \popt \space procedure has prior mean power almost as high as the constrained conditional \tk \space procedure, but the second row indicates that it also selects the more robust design (i) far more often (usually when the posterior probability of $H_1$ is low). Thus, the conditional \popt \space procedure offers a compromise between power for distinguishing the hypotheses of interest and robustness.


\section{Links between Test and Estimation Information}
\label{sec:theory}

\subsection{Fraction of observed test information}
\label{sec:test_frac}

\cite{frac_info} proposed several measures of the fraction of observed test information to guide data collection decisions in genetic studies (see Section \ref{sec:uses}).   We provide the general mathematical form in Definition \ref{def:frac_dinfo} (below) because the fraction of observed test information is important in sequential design and for establishing theoretical connections between test and estimation information.

\begin{defi}
\label{def:frac_dinfo}
The fraction of observed test information provided by the first part of the composite experiment $\ept=(\ept_1,\ept_2)$ for comparing the hypotheses $H_0: \para\in \psupnull$ and $H_1: \para \in \psupalt$, for a given evidence function $\declr$  and a proper prior $\prior$, is defined as
\begin{align}
\fracinfo(\ept_2|\ept_1;H_0,H_1,\prior,\datx_1) =\frac{\dinfo(\ept_1;\datx_1)}{\dinfo(\ept_1;\datx_1)+\dinfo(\ept_2|\ept_1;\datx_1)},\label{eqn:frac_dinfo}
\end{align}
where $\datx_1$ is the observed outcome of $\ept_1$, and $\psupnull \cap \psupalt = \emptyset$.
\end{defi}

If $\dinfo(\ept_1;\datx_1) \geq 0$,  then it follows that (\ref{eqn:frac_dinfo}) is  between 0 and 1. In practice, if $\fracinfo(\ept;\datx_1)$ is close to one then we may decide not to perform $\ept_2$, particularly if it is expensive. The canonical example sets $\declr(\zlr)=\log(\zlr)$ and thus takes the ratio of the observed data log Bayes factor and the expected complete data log Bayes factor. Similarly, in the frequentist case, \cite{frac_info} suggested the measure
\begin{align}
\label{eqn:ri1}
\ri_1=\frac{\log \lr(\para\ob,\para_0|x\ob)}{\Exmis[\log \lr(\para\ob,\para_0|\ranx\ob,\ranx\mis)|\para\ob,\datx\ob]},
\end{align}
where $\para\ob$ is the MLE of $\para$ based on $\datx\ob$. 

The decision whether to collect more data depends on which hypothesis is true, because if the observed data supports the false hypothesis then our need for additional data is greater. Thus, it is unsurprising that there is no general coherence identity for the fraction of observed test information. In practice, we suggest using $\fracinfo(\ept_2|\ept_1;H_0,H_1,\prior,\datx_1)$ if $\zlr(\datx_1) \leq 1$ and $\fracinfo(\ept_2|\ept_1;H_1,H_0,\prior,\datx_1)$ otherwise. The resulting measure has a similar interpretation as (\ref{eqn:frac_dinfo}) but takes account of which hypothesis is more likely, and is always between $0$ and $1$. 
In the special case where $\declr(1)=0$ (and (\ref{eqn:observe_coherence}) is satisfied), we have $\fracinfo(\ept_2|\ept_1;H_0,H_1,\prior,\datx_1)=\fracinfo(\ept_2|\ept_1;H_1,H_0,\prior,\datx_1)$, but this is not a coherence identity since the corresponding observed test information is negative on one side of the equality and positive on the other.

\subsection{Connections between estimation and test information}
\label{sec:asymptotics}

\cite{meng_dyk1996} showed that the {\it relative augmentation function}
\begin{align}
\ri(\para)=\frac{\log \lr(\para\ob,\para|x\ob)}{\Excom[\log \lr(\para\ob,\para|\ranx\co)|\para\ob,x\ob]}
\end{align}
converges to the fraction of observed Fisher information
\begin{align}
\ri_E = \frac{I\ob}{I\ob+I\mis}\label{eqn:re}
\end{align}
as $|\para -\para\ob| \rightarrow  0$. Here, $I\ob$ is the usual observed Fisher information,
and $I\mis$ is the missing Fisher information given by
\begin{align}
I\mis = \Excom\left[\left.\left.-\frac{\partial^2\log f(\ranx\co|\datx\ob,\para)}{\partial \para^2}\right|\datx\ob,\para\right]\right|_{\para=\para\ob}. \label{eqn:missing_info}
\end{align}
As \cite{frac_info} mentioned, replacing $\para$ with $\para_0$ gives us the same limit for the measure $\ri_1$ in (\ref{eqn:ri1}). This result is intuitive in that we might expect test information to coincide with estimation information when the two hypotheses are both very close to $\para\ob$. The following theorem generalizes the equivalence, and its proof is given in Appendix \ref{app:proof}.

\begin{thm}\label{theorem} Let the hypotheses be $H_0:\para=\para_0$ and $H_1:\para=\para_1$, and suppose that the derivatives of the evidence function $\declr$ exist at 1. Then, for univariate $\para$ and $\para_1=\para\ob$, we have
\begin{align}
\fracinfo(\ept_2|\ept_1;H_0,H_1,\datx\ob) = \frac{\declr'(1)I\ob}{\declr'(1)I\ob-\declr''(1)I\mis} + O_p(|\para_0 - \para\ob|), \label{eqn:fractest_limit}
\end{align}
under the uniform integrability condition given in the proof in the Appendix.
\end{thm}

It is possible to extend Theorem~\ref{theorem} to avoid the univariate condition and sharp hypotheses (by using priors that converge to delta functions), but the current form suffices to illustrate the connection between test and Fisher information. The theorem tells us that if $\declr'(1)=-\declr''(1)$, then  $\fracinfo(\ept_2|\ept_1;H_0,H_1,\datx\ob)$ will exactly correspond to $\ri_E$ as $|\para_0 -\para\ob|\rightarrow 0$. Otherwise, the {\it relative conversion number} $\cn = -\declr''(1)/\declr'(1)$ indicates how much of the missing data estimation information is converted to test information in the limit, relative to the conversion of observed estimation information. 
For example, under the \postv evidence function (\ref{eqn:postef_h0}) we have $\cn=2\pi_0$, and therefore the stronger our initial bias in favor of the null, the greater the importance of the missing data estimation information, relative to the observed estimation information.

This makes sense because Fisher information measures our ability to estimate the true parameter, and the value of successful estimation for testing depends on the strength of our prior separation of the hypotheses.  If $\pi_0=0.5$, then all estimation information will be helpful because the prior does not separate the hypotheses, hence $\cn=1$  and the fraction of observed test and estimation information coincide. When $\pi_0$ is close to $0$, the posterior probability of $H_1$ (the hypothesis assumed true by Definition \ref{def:frac_dinfo}) will be close to one, even though the observed data provides no evidence. Thus, we have $\cn\approx0$, because there is little to be gained by collecting more data. When $\pi_0$ is close to $1$, the prior is in conflict with our assumption that $H_1$ is true, and therefore estimation information from new data not only has the potential to distinguish the hypotheses, but also to overcome misleading information from the prior. Hence, when $\pi_0$ is close to $1$, the value of estimation information is inflated in the testing context, and we expect $\cn > 1$. In the current example $\cn=2\pi_0$,  meaning that indeed $\cn \approx 2 > 1$, as anticipated, when $\pi_0$ is close to $1$. 

The relative conversion number  has similar interpretations for other evidence functions. In each case it indicates the relative worth of the missing data estimation information for testing, when there is no evidence in the observed data. Thus, $\cn$ provides a characterization of the general approach to testing implied by the evidence function, i.e., whether we would be likely to collect additional data if the observed data does not separate the hypotheses. 

\section{Follow-up Observations of Astronomical Lightcurves}
\label{sec:application}

\begin{figure}[t]
\includegraphics[width=1\textwidth,trim = 0mm 0mm 0mm 5mm]{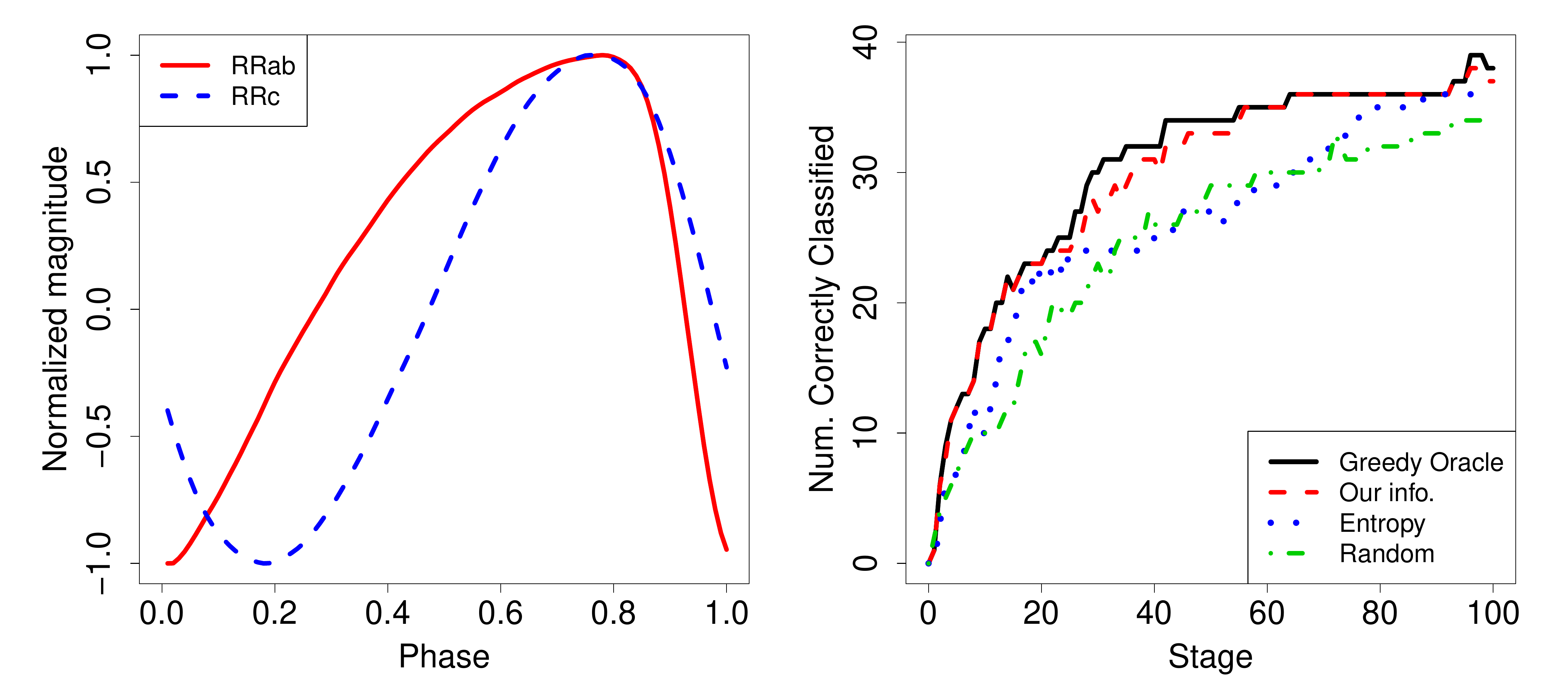}
\caption{Left: lightcurve templates for RRab (solid line) and RRc (dashed line) stars. Right: number of  lightcurves that are correctly classified at each stage of follow-up observations. The number is out of the 61 lightcurves that were initially misclassified based on posterior probability. The four lines correspond to the four different methods for scheduling the follow-up observations discussed in the text. \label{fig:lcresults}}
  \end{figure}

We now demonstrate how our methodology can be used to schedule follow-up telescope observations in order to better distinguish two subclasses of RR Lyrae stars. RR Lyrae stars pulsate, which leads to periodic changes in their brightness, and they are important objects because they allow astronomers to calculate distances to other galaxies and within our own galaxy. 
Our dataset consists of 200 RRab and 200 RRc lightcurves from the Catalina Real Time Transient Survey (CRTS), see \citet{drake2009} and \citet{drake2014}. A lightcurve is a time series recording the brightness of a star over time, and different characteristic lightcurve shapes are associated with each type of star. The left panel of Figure \ref{fig:lcresults} shows templates we constructed for RRab and RRc lightcurves from a large training dataset again from CRTS. The templates are plotted in phase space because the lightcurves of RR Lyrae are periodic (up to observation noise). 
The full CRTS catalog is available at \url{http://crts.caltech.edu/}.

Our dataset consist of 400 lightcurves, each with between 62 and 499 observations. Through initial screenings, we can assume they are either RRab or RRc. Given a statistical model for the two
types, we can compute the posterior probability that a lightcurve is RRab.  Sometimes the posterior probability will be close to 0.5 or favor the wrong lightcurve type (here we know the true types for practical purposes because experts have  inspected each lightcurve individually and taken additional information into account).  It is therefore of interest to know what follow-up observation times would allow us to better identify the true class. In practice, telescope time is limited and so we suppose that for each lightcurve there are 3 possible follow-up observation times to choose from. We will choose using our probability based information measure, which uses the evidence function (\ref{eqn:postef_h0}), and make comparisons with several other methods.

First we describe the lightcurve model. The templates in the left panel Figure \ref{fig:lcresults} are mean lightcurve shapes and therefore we expect deviation from these shapes in actual observations due to both variability in the underlying lightcurve shapes and observation noise. To account for these sources of variability we model each lightcurve as a Gaussian process (GP) with the mean being the appropriate template and a nugget or additional variance term to account for observation noise. Specifically, let $x_k = (x_{k1},\dots,x_{kn_k})$ be the magnitude (brightness) observations for source $k$, with corresponding observation phases  $t_k = (t_{k1},\dots,t_{kn_k})$, for $k=1,\ldots,400$. We use observation phases rather than times because for the lightcurve in our dataset it is possible to estimate the period accurately so we treat the period as known. If $\rho$ is the lightcurve period, then the phase at time $t$ is given by $(t \text{ mod } \rho)/\rho$. (Note that sometimes period  information alone can distinguish between types, but this is not always the case for RRab and RRc stars, which motivates our investigation.)
Under class $C=c$, we model the $k^{\text{th}}$ lightcurve observations by the GP
\begin{align*}
X_k| t_k, C=c \sim N(\alpha_{ck}+s_{ck}\mu_c(t_k), \Sigma_{ck}),
\end{align*}
where $\mu_c(t_k)$ is a vector whose entries are the values of the template for class $c$ at the phases given by $t_k$, and $\Sigma_{ck}$ is an $n_k \times n_k$ covariance matrix, for $k=1,\dots,400$. The parameters $\alpha_{ck}$ and $s_{ck}$ center and scale the lighcurve template. Let $S_b=[(b-1)/10,b/10)$, for $b=1,\dots,10$, be disjoint phase bins on which we assume fixed nugget terms $\tau_{ckb}^2$, for $k=1,\dots,400$. We then decompose $\Sigma_k$ as $V_{ck} + D_{ck}$, where
\begin{align}
\left(V_{ck}\right)_{ij}  = \sigma_{ck}^2\exp\left(-\frac{(t_{ki}-t_{kj})^2}{2\lambda_{ck}^2}\right),\label{eqn:kernel}
\end{align}
and $D_{ck}$ is a diagonal matrix with $i^{\text{th}}$ diagonal entry given by $\tau_{ckb}^2$, where $b$ is such that $t_{ki}\in S_b$. For simplicity and computational reasons, for each lightcurve we find an initial estimate of the alignment and noise parameters $\phi_{ck}=(\alpha_{ck}, s_{ck}, \tau_{ck1},\dots,\tau_{ck10})$, and then treat these parameters as known. However, the GP parameters $\beta_{ck} =(\sigma_{ck},\lambda_{ck})$ are nuisance parameters to be fitted. We use a Laplace approximation for the observed Bayes factor $\bfac(x_k|H_0,H_1)$ (see (\ref{eqn:nuisance_bf})), where $H_0:c=1 ({\rm RRab}) $ and $H_1:c=2 ({\rm RRc})$.

At each stage $s$, for each lightcurve $k$, we randomly generate 3 phases $t^{\text{new}}_{sk1}, t^{\text{new}}_{sk2}, t^{\text{new}}_{sk3}$, for follow-up observations and select one of them using four different methods: (i) an oracle method, (ii) the conditional test information under the probability based evidence function (\ref{eqn:postef_h0}) (with $\pi_0=\pi_1=0.5$), (iii) the entropy based method of \citet{boxhill1967}, and (iv) random selection. The oracle method selects the observation phase that results in the highest expected posterior probability for the true class (which is otherwise treated as unknown).
We use Monte Carlo approximations to perform the averaging needed for methods (i)-(iii). After choosing a future phase under each method, the follow-up observations at the new phases are simulated under the model described above, with $c$ set to the true class. We repeat this process for stages $s=1,\dots,100$. At each stage, the simulations of the new observations are performed conditioning on the initial lightcurve and the data points collected up to that stage under the oracle method.

For the initial lightcurves, 61 of the 400 stars were misclassified when assigning them to the class with the highest posterior probability. In this example the stars that were initially correctly classified were still correctly classified when follow-up observations were collected, so for the purposes of evaluation it makes sense to focus on only the initially misclassified ones. The right panel of Figure \ref{fig:lcresults} shows the number of the 61 misclassified lightcurves that were correctly classified at each of the 100 follow-up stages under each of the four scheduling methods. The random method performs worst because the dash-dot line is below the other lines at most stages, meaning that to achieve a given number of correct classifications the random method generally requires more follow-up observations. The oracle method (solid line) unsurprisingly performs the best, and our test information approach (dashed line) performs very similarly. In fact, our probability based test information should theoretically match the oracle method, due to the coherence identity (\ref{eqn:con_coherence}), which ensures we do not need to know the true class to choose the optimal strategy. The discrepancy seen in Figure \ref{fig:lcresults} is due only to numerical inaccuracies, i.e., we use a Laplace approximation to our test information under the null and under the alternative, but the two do not exactly match.  The entropy method of \citet{boxhill1967} (dotted line) performs better than the random approach but not as well as our test information. In Section \ref{sec:discussion} and Appendix \ref{app:entropy}, we provide conceptual reasons to expect our information to offer improvement over the entropy method.

\section{Discussion and Further Work}
\label{sec:discussion}

\subsection{Philosophical considerations}

\citet{boxhill1967} proposed maximizing the change in the expected entropy of the posterior distribution on the model space. In the current context, their criterion is given by
\begin{align}
\phi_{BH}(t) = -\sum_{i=0}^1 \pi_i \log \pi_i + \Ex\left.\left[P(H_i|X,t)\sum_{i=0}^1P(H_i|X,t)\right|t\right], \label{eqn:bh_crit}
\end{align}
where  $t$ is the design point, and $P(H_i|X,t)$ is the posterior probability of $H_i$, given the data $X$ (collected at the design point $t$). Importantly, the expectation in (\ref{eqn:bh_crit}) is taken  with respect to $f(\cdot|t)=\pi_0 f(\cdot|H_0,t) + \pi_1 f(\cdot|H_1,t)$, where $f(\cdot|H_i,t)$ is the conditional probability density function of $X$, given the hypothesis $H_i$ and the design point $t$, for $i=1,2$. A limitation of this approach is that there is an inconsistency in taking the expectation with respect to $f(\cdot|t)$, which averages over the two hypotheses, because in most scientific scenarios it is assumed that only one of the hypotheses is true. From a Bayesian perspective, which assumes a distribution on the unknown true hypothesis, maximizing (\ref{eqn:bh_crit}) is the right approach. However, the hypothesis testing questions posed in science are not Bayesian, and therefore optimizing (\ref{eqn:bh_crit}) can lead to sub-optimal decisions. For example, it is possible to construct a situation where  (\ref{eqn:bh_crit}) favors one design but another design is uniformly better in terms of the probability of coming to the right conclusion (even if the conclusion is to be drawn based on the posterior probability of each hypothesis). See Appendix~\ref{app:entropy} for an example, which is somewhat contrived, but nevertheless proves the existence of this phenomenon.

Our information criteria avoid the above issue by exploiting the fundamental coherence identities given in (\ref{eqn:info_link}) and (\ref{eqn:con_coherence}), which allow us to average over the future data while conditioning on only one hypothesis and then make use of the symmetry between the hypotheses.
In the application in Section \ref{sec:application}, our probability based information criterion performs somewhat better than the  \citet{boxhill1967} criterion, and chose the same design points as the oracle method more often. It is possible that this superior performance is due to the use of the coherence identities to avoid the issues seen in Appendix \ref{app:entropy}, though more investigation is clearly needed.

\subsection{Additional benefit}
\label{sec:power}


Our general framework can also help to rule out undesirable  measures of test information. For example, the variance of the log Bayes factor is essentially a measure that \citet{frac_info} rejected after some trial and error, but our framework will reject it immediately because the corresponding evidence function $\declr(\zlr)=-(\log z - c)^2$ is not concave. Perhaps the most notable quantity ruled out is power, which has no easy quantification in the presence of composite hypotheses or nuisance parameters. 
The second major problem is that there is not an intuitive measure of {\it observed} power, and it is therefore 
difficult to see how  sequential design decisions (e.g., stopping rules) can be based on power. A further fundamental difficulty is that power does not have the maximal information interpretation discussed in Section \ref{sec:expected_test} because it incorporates an investigator-specific critical region.

\subsection{Future work}

A natural future step is to investigate how test and estimation information measures can be combined to find designs that are good for both testing and estimation. Some work has been done along these lines by \citet{borth1975} in the special case of the entropy approach taken by \citet{lindley1956} and \citet{boxhill1967}. However, in general, test and estimation information are not related simply, and therefore trying to directly find designs that are good for both testing and estimation may not be an effective strategy. Instead, we can divide up the design points and construct two designs, one that is good for testing and one that is good for estimation. The overall design, composed of all the design points, should then have reasonable properties for both tasks. Future work will explore this approach and investigate methods for setting the proportion of the design points to be allocated to each task. 

Another key direction is to deal with more than two hypotheses, a scenario which often arises in classification contexts. In preliminary investigations, we have extended many of the concepts in this paper to multiple hypotheses by using one-versus-all comparisons, but more work is needed to demonstrate the utility of our framework in this more general setting.

\bibliographystyle{chicago}
\bibliography{references.bbl}

\appendix

\section{Proof of Theorem 1}\label{app:proof}

A Taylor expansion of $\declr(1)-\declr\left(f(\datx|\para_0)/f(\datx|\para_1)\right)$ in $\para_0$ gives
\begin{align}
-\left( \delta\frac{f'(\datx|\para_1)}{f(\datx|\para_1)}  + \frac{\delta^2}{2}\frac{f''(\datx|\para_1)}{f(\datx|\para_1)}\right)\declr'(1)-\frac{\delta^2}{2}
\left(\frac{f'(\datx|\para_1)}{f(\datx|\para_1)}\right)^2\declr''(1)+ R_3(\datx),\label{eqn:expandv}
\end{align}
where $\delta = (\para_0-\para_1)$, and $R_3(\datx)= R_3(\datx;\para_0,\para_1)$ is the standard Taylor expansion remainder term. If we set $\para_1=\para\ob$, then $\dinfo(\ept_1;\datx\ob)$  becomes
\begin{align}
\frac{1}{2}(\para_0-\para\ob)^2I\ob\declr'(1) + O\left((\para_0-\para\ob)^3\right).\label{eqn:app_obs_info}
\end{align}
Next, for a sequence $\{\para^{(m)}\}_{m=1}^\infty$ such that $|\para^{(m)}-\para_1| \leq \frac{1}{m}$, we assume uniform integrability for the sequence $\{R_3(\{\datx\ob,\ranx\mis\};\para^{(m)},\para_1)\}_{m=1}^\infty$. Then, inserting $(\ranx\mis,\datx\ob)$ for $\datx$ in (\ref{eqn:expandv}), setting $\para_1=\para\ob$, and taking an expectation with respect to $f(\ranx\mis|\datx\ob,\para\ob)$, we obtain
\begin{align}
\frac{1}{2}(\para_0-\para\ob)^2\left(I\ob\declr'(1)-I\mis\declr''(1)\right)+ O_p\left((\para_0-\para\ob)^3\right).\label{eqn:complete_info}
\end{align}
The result (\ref{eqn:fractest_limit}) then follows after taking the ratio of (\ref{eqn:app_obs_info}) and (\ref{eqn:complete_info}) and letting $\para\ob\rightarrow \para_0$.


\section{Entropy criterion example} \label{app:entropy}

Suppose there are two choices for a design point $t^*$, say $t_1$ and $t_2$.
Here we show that it is possible to construct an example such that
 \begin{align}
  \phi_{BH}(t_1)&< \phi_{BH}(t_2) \label{eqn:bh1}\\
    P(P(H_T|X,t_1)>0.5 | t_1)& \geq  P(P(H_T|X,t_2)>0.5| t_2) \label{eqn:bh2}\\
 \Ex[P(H_T|X,t_1) | t_1]& >   \Ex[P(H_T|X,t_2)| t_2]\label{eqn:bh3}\\
 \phi_P(t_1)&> \phi_P(t_2)\label{eqn:bh4}\\
   P(P(H_T|X,t_1)>0.5 | t_1)& >  P(P(H_T|X,t_2)>0.5| t_2) \hspace{0.1cm}\text{for some $H_T$},\label{eqn:bh5}
 \end{align}
 where $\phi_{BH}$ is the entropy based criterion of \citet{boxhill1967} given in (\ref{eqn:bh_crit}), $X$ is the data to be sampled under $t^*$, $H_T$  is the true hypothesis,  and the (\ref{eqn:bh1})-(\ref{eqn:bh4}) hold regardless of whether $H_T=H_0$ or $H_T=H_1$.
 In other words, there are situations where the criterion $\phi_{BH}$ favors $t_2$ even though the  posterior probability of $H_T$ is more likely to be greater than 0.5 under $t_1$ than under $t_2$ for at least one choice of $H_T$, and at least as likely for the other choice of $H_T$ (see (\ref{eqn:bh5}) and (\ref{eqn:bh2})). Furthermore, (\ref{eqn:bh3}) and (\ref{eqn:bh4}) convey that the expected posterior probability for the true hypothesis and our probability based information criterion, respectively, are also larger under $t_1$, which lends further support to favoring $t_1$.

{\it Example.} Suppose that  $H_0:X \sim f_0 \equiv \beta_1\text{Uniform}(0,1) + (1-\beta_1)\text{Uniform}(1,2)$, regardless of $t^*$. Under $H_1$, let the distribution of $X$ depend on $t^*$ as follows:
  \begin{align*}
 X\sim \left\{\begin{array}{lc}\alpha f_0  + (1-\alpha)\text{Uniform}(3,4) & \mbox{if }t^*=t_1,\\
  \beta_2\text{Uniform}(0,1) + (1-\beta_2)\text{Uniform}(1,2) & \mbox{if }t^*=t_2.
  \end{array}\right.
  \end{align*}
Setting $\pi_0=0.999, \pi_1=0.001, \alpha=0.99$, $\beta_1=0.1$, and $\beta_2=0.9$, the inequalities (\ref{eqn:bh1})-(\ref{eqn:bh4}) hold for $H_T=H_i$, for $i=0,1$, and $(\ref{eqn:bh5})$ holds for $H_T=H_1$. The key to this example is that  $\phi_{BH}(t_1)$ is small despite the fact that the only possibility for $P(H_1|X,t_1)>0.5$  occurs when $H_T=H_1$ and $t^*=t_1$. It is easy to construct a similar example using Normal distributions.

\end{document}